\documentclass[11pt]{article}
\usepackage{setspace}
\usepackage[pdftex]{graphicx}
\usepackage{amssymb,amsmath}
\newdimen\IEEEelabelindent
\IEEEelabelindent \parindent

\def\beq{\begin{equation}}
\def\eeq{\end{equation}}
\def\bqn{\begin{eqnarray}}
\def\eqn{\end{eqnarray}}

\def \ee {\begin{equation}}
\def \eee {\end{equation}}
\def \eqe {\begin{eqnarray}}
\def \eqee {\end{eqnarray}}

\newtheorem{thm}{Theorem}[section]
\newtheorem{positive}[thm]{Lemma}
\newtheorem{kernel}[thm]{Lemma}
\newtheorem{isolated}[thm]{Theorem}
\newtheorem{stable}[thm]{Theorem}
\newtheorem{lawrence}[thm]{Theorem}
\newtheorem{pairing}[thm]{Lemma}
\title{Tensegrity Models and Shape Control of Vehicle Formations}
\author{Benjamin Nabet$^*$ and Naomi Ehrich Leonard\thanks{B.~Nabet and N.E.~Leonard are with the Department of Mechanical and Aerospace Engineering, Princeton University, Princeton, NJ, 08544 USA {\tt\small \{bnabet,naomi\}@princeton.edu.}  This work is supported in part by NSF grant  CMS-0625259 and ONR grants
N00014--02--1--0826 and N00014--04--1--0534. }}
\begin{document}
\maketitle

\begin{abstract}
Using dynamic models of tensegrity structures, we derive provable, distributed control laws for stabilizing and changing the shape of a formation of vehicles in the plane.   Tensegrity models define the desired, controlled, multi-vehicle system dynamics, where each node in the tensegrity structure maps to a vehicle and each interconnecting strut or cable in the structure maps to a virtual interconnection between vehicles.  Our method provides a smooth map from any desired planar formation shape to a planar tensegrity structure.  The stabilizing vehicle formation shape control laws are then given by the forces between nodes in the corresponding  tensegrity model.   The smooth map makes possible provably well behaved changes of formation shape over a prescribed time interval.  A designed path in shape space is mapped to a path in the parametrized space of  tensegrity structures and the vehicle formation tracks this path with forces derived from the time-varying tensegrity model.  By means of examples, we illustrate the influence of design parameters on performance measures.

\end{abstract}

\section{Introduction}

\qquad Recent efforts in coordinated control have focused on
reconfigurable mobile sensor networks with application to collective sensing and monitoring on land, in the sea, in the air and in space.  Each sensor platform in the network is a vehicle; therefore, the spatial distribution of the vehicles determines the geometry and resolution of the sensor array.   Control over the geometry and resolution of the vehicle formation, also referred to as the ``shape'' of the formation,  can provide important advantages to performance and efficiency of data gathering and processing.   For example, in searching and tracking tasks, it can be critical  to have the network estimate gradients and possibly higher-order derivatives from noisy measurements of the sampled field.   To minimize the error in these estimates, the shape  of the formation should adapt to changes in the environment and the motion of the network.    Zhang and Leonard presented an algorithm for level set tracking where the shape of the group was determined so as to minimize the least mean square error in gradient estimates of a scalar field \cite{swarm}.    More generally, shape control plays an important role in multi-scale sensing tasks where changes are desired that require shifts in resolution for some or all of the network.   Shape control can be significant in other vehicle network tasks, for example, when vehicles coordinate their activity to escort, carry or otherwise interact with objects in their environment.

The shape of a vehicle formation depends only on relative measurements among individuals, i.e., it describes the way the individual vehicles are arranged relative to one another rather than where the group is or how it is oriented.   Accordingly, shape control can naturally be distributed.    In this paper, we propose a methodology that systematizes the
design and analysis of distributed shape control of a group of vehicles in the plane.

Our approach is to synthesize and analyze shape dynamics of a group using models of {\em tensegrity structures}, which are
spatial networks of interconnected struts, cables and rods
\cite{snelson,fuller,intro,connelly1}. The goal is to design
control laws that drive vehicle formations into shapes with forces
that can be represented as those internal to tensegrity
structures. We extend the concept to control the formation
along a smooth path in shape space consisting of tensegrity structures;
this allows for stable, smooth reconfiguration of the group
shape.

The artist Kenneth Snelson \cite{snelson} built the first tensegrity structure, and Buckminster Fuller \cite{fuller} coined the term ``tensegrity" by combining the words ``tension" and ``integrity".  Tensegrity structures have since been extensively studied in the engineering, mathematics and even the biology literatures; see \cite{Sultan06} and references therein for a history and survey of tensegrity research.   In the mathematics literature, Connelly, Whiteley and others have studied and proved conditions for stability and rigidity properties of tensegrities (sometimes referred to in the mathematical formulation as tensegrity frameworks) \cite{connelly1,connelly2,connelly3,connelly4}.    These properties help motivate our approach since  modeling a multi-agent network as a tensegrity suggests that the vehicle network will inherit the same useful properties as the tensegrity.  In the engineering literature, researchers have studied both statics and dynamics of tensegrity structures, e.g., in \cite{Motro86,Sultan02}, and have explored numerous possibilities for tensegrity structures including both passively and actively controlled structures, e.g., in \cite{Sultan97,Motro98}.
Skelton et al  \cite{intro} define the ``small control energy principle" arguing that the shape of a tensegrity structure can be changed with little control energy; their method involves changing the physical parameters of the structure so that with little effort the shape of the structure changes  to the new corresponding equilibrium.   Our approach to shape change similarly makes the desired shape a new equilibrium for the system.   In \cite{SultanIJSS02} the authors study symmetrical motions for reconfiguration of a class of tensegrity structures composed of six bars, a rigid top and base and eighteen tendons.    

The tensegrity structures we consider are formed by a combination of \emph{struts} (connecting elements in compression) and \emph{cables} (connecting elements in tension), which we classify together more generally as \emph{edges}. These edges meet at \emph{nodes}. A significant challenge in the tensegrity literature has been solving the \emph{form finding} problem: determining the shape of the tensegrity, given a number of nodes and edges \cite{formfinding1}.  For the ``reverse engineering" problem, i.e., determining a model and a set of edges, given a desired shape, Connelly has proven a result that provides a means to systematically design stable planar tensegrities in the shape of any strictly convex polygon \cite{connelly1}. The method developed by Connelly is inductive (i.e., a tensegrity with $N$ nodes is constructed from a tensegrity with $N-1$ nodes) and yields structures that almost never require all-to-all coupling topologies.   In the present paper, we develop an alternative systematic method to define a tensegrity model that realizes any given arbitrary planar shape that is convex or non-convex, at the cost, however, of often requiring connecting topologies closer to all-to-all coupling.   Our method also has the advantage of providing a mapping from shape to tensegrity that is smooth with respect to shape parameters.   

Using an energy approach, Connelly developed a model of the forces and formalized the notion of stability of a tensegrity \cite{connelly1,connelly2,connelly3}.  In his model it is assumed that cables do not increase in length and struts do not decrease in length.  The tensegrities we design to model a vehicle network are virtual since the vehicles are not physically connected. Accordingly,  these constraints that are reasonable in the physical setting have to be relaxed in the virtual setting, since restrictions on cable and strut length changes  cannot be imposed as a constraint on the distances between pairs of vehicles that are not physically connected.  In the absence of the cable and strut length change constraints, the model from \cite{connelly1}, which uses linear springs to represent cables and struts, allows for arbitrary stretching and shrinking of a tensegrity in the plane.   We therefore augment the model with the relaxed constraints in order to isolate planar tensegrities. We then present a systematic and smooth method to compute the parameters of the augmented model that realizes a tensegrity structure for any arbitrary desired shape in the plane.

To make our mapping from shape to tensegrity relevant for vehicle network control, we define the following mapping between the vehicle network and a tensegrity structure: each node of the
tensegrity structure is identified with one vehicle of the network
and the edges of the structure correspond to communication and control force
directions  between the vehicles. If an edge is a cable,
the force is attractive; if the edge is a strut then the force is
repulsive. The magnitude of the force depends on the tensegrity
parameters and the relative distance between the
vehicles associated with the edge. The dynamics of the tensegrity nodes are derived assuming a point-mass double integrator model for each networked vehicle. The point-mass model
may appear somewhat simplistic since it seems to ignore the
challenges of controlling the detailed dynamics of each of the
individual vehicles. However, for practical
implementations such as the coordinated control of a network of autonomous underwater gliders for adaptive ocean sampling, as described in \cite{Leonard07},  it is useful to decouple
the design of the coordinating control strategy from the
lower-level, individual-based tracking control. 

Using the mapping we derive between shape and tensegrity and the mapping we define between tensegrity and vehicle network, we prove that the tensegrity model provides control forces that stabilize a network of $N$ vehicles to a desired shape.   We extend this approach to stable, smooth vehicle network shape change control; the method finds and tracks a smooth path in the space of tensegrities (equivalently the space of shapes) that connects prescribed initial and final shapes.

In Section~\ref{sec:models} we present models of tensegrity including our augmented model to isolate tensegrity shapes without constraining edge lengths.   We then derive a smooth map from arbitrary planar shape to planar tensegrity.   In Section~\ref{sec:stability} we prove stability of the desired shape for the dynamics of the vehicle network represented as a tensegrity. We prove and illustrate our method for controlling change of shape of the network in Section~\ref{shape}.  We make final remarks in Section~\ref{sec:finalremarks}.

\section{Tensegrity Models and Mapping Arbitrary Shape to Tensegrity}
\label{sec:models}
\subsection{Shape equilibrium conditions for model with edge length constraints}

\qquad We investigate in this section the modeling of forces induced by the two types of edges of a tensegrity structure as developed by Connelly in \cite{connelly1}. With this model and physical constraints imposed on the lengths of struts and cables, Connelly has proven a theorem that provides a
means to systematically design stable planar tensegrities that correspond to strictly convex polygons \cite{connelly1}.  As a first step to developing our systematic and smooth method for designing planar tensegrities to realize arbitrary shapes, we first examine the equilibrium conditions for  Connelly's dynamic model.  Using  algebraic graph theory, we formally derive the conditions on allocation of cables and struts and choice of model parameters that must be satisfied for a given shape to be an equilibrium of the dynamic model.  We show that in the absence of the cable and strut length constraints, this simple model, which uses linear springs to represent cables and struts, has a continuum of equilibria allowing arbitrary stretching of a tensegrity in the plane. Hence in our virtual setting, where physical constraints on the edges cannot be imposed, the linear model does not produce an isolated
equilibrium as we desire; in particular, the model cannot be used to control the scale of the group shape.  In Section~\ref{nonlinear} we augment the model to isolate shapes of prescribed scale.

The edges of a tensegrity structure are modeled as linear springs
with zero rest length. Cables have a positive spring constant
while struts have a negative one \cite{connelly1,connelly2}. Hence
for two nodes $i, j$ we have
\begin{equation}  \label{force:model1}
\vec{f}_{i\rightarrow j} = \omega_{ij}(\vec{q}_i-\vec{q}_j)=-
\vec{f}_{j\rightarrow i},
\end{equation}
where $\vec f_{i\rightarrow j} \in \mathbb{R}^2$ is the force
applied to node $j$ as a result of the presence of node $i$,
$\vec{q}_i = (x_i,y_i) \in \mathbb{R}^2$ is the position vector of
node $i$, and $\omega_{ij}$ is the spring constant of the edge
$ij$. The spring constant $\omega_{ij}$ is negative if $ij$ is a
strut, positive if $ij$ is a cable, and zero if there is no
connection between the nodes $i$ and $j$; $\omega_{ij}$ is called the
\emph{stress} of the edge $ij$. For a tensegrity structure with
$N$ nodes, its weighted interconnection topology can be described by the vector $\boldsymbol{\omega}\in\mathbb{R}^{\frac{N(N-1)}{2}}$, where components consist of the
stresses $\omega_{ij}$. If the vector $\boldsymbol{\omega}$ has no zero elements, then the associated interconnection graph is
complete, i.e., there is an edge between every pair of nodes. A given shape can be realized with different sets of edges, notably those with incomplete graphs.

Shape refers to the way the nodes are arranged relative
to one another rather than where the structure is or how it is
oriented. Thus, a given shape can be associated with different sets of absolute position vectors. Let $\mathbf{x}=(x_1,\ldots,x_N)^T$ and
$\mathbf{y}=(y_1,\ldots,y_N)^T$, then a {\em placement}, defined as $\mathbf{q}=\begin{pmatrix}
\mathbf{x} \cr \mathbf{y}\end{pmatrix}\in \mathbb{R}^{2N}$, characterizes the tensegrity in absolute coordinates (for convenience we will write  $\mathbf{q} = (\mathbf{x},\mathbf{y}))$.   Let $\vec{q}_c = (1/N)\sum\limits_{i=1}^N \vec{q}_i$ be the centroid of the tensegrity.  For two placements $\mathbf{q}^1 \in \mathbb{R}^{2N}$ and $\mathbf{q}^2 \in \mathbb{R}^{2N}$, we define the equivalence relation $\mathcal{R}$ as
\begin{equation*}
\mathbf{q}^1\mathcal{R}\mathbf{q}^2\Longleftrightarrow \exists (R, \vec{t})\in
SE(2) \;\; \mbox{and}\;\;
\mathbf{q}^2 = (\mathbf{q}^1 \; \mbox{after rigid rotation about} \; \vec{q}_c^1 \; \mbox{by} \; R \; \mbox{and translation by} \; \vec{t}).
\end{equation*}
An element of $SE(2)$ is a rigid motion and so any two placements in the same equivalence class  $\mathcal{R}$ have the same shape.   We can therefore identify a given shape with the equivalence class $[\mathbf{q}^1] = \{ \mathbf{q} \in \mathbb{R}^{2N} \; | \; \mathbf{q}^1 {\cal R} \mathbf{q}\}$, where $\mathbf{q}^1$ is a representative placement with the given shape.   The class $[\mathbf{q}^1]$ can be identified with  $SE(2)$.  We note that the ordering of the nodes in the placement matters to shape classification since two placements $\mathbf{q}^1$ and $\mathbf{q}^2$ representing the same geometric shape but with nodes permuted  will not be in the same equivalence class.

Using the force model (\ref{force:model1}), we derive the
equations of motion for each node of a tensegrity with placement $\mathbf{q}$.   The potential energy of a
tensegrity structure, with the forces induced by its edges as modeled
by (\ref{force:model1}), is
\begin{equation} \label{potential1}
V(\mathbf{q})=\frac{1}{2}\sum_{i=1}^{N}\sum_{j=i+1}^{N}\omega_{ij}\|\vec{q}_j-\vec{q}_i\|^{2}.
\end{equation}
In the sequel we write  $\sum\limits_{i<j}$ to represent
$\sum\limits_{i=1}^{N}\sum\limits_{j=i+1}^{N}$. The potential (\ref{potential1})
increases as we stretch the cables and shrink the struts. We
introduce in the system a linear damping force of the form
$-\nu\dot{\vec{q}}_i,$ where $\nu>0$ is a damping coefficient.
The Euler-Lagrange equations of motion for the
nodes of the tensegrity structure using Cartesian coordinates in
the plane are
\begin{equation}\label{eomfirstorder}
\left\{ \begin{aligned}
\dot{x}_i&=p_i^x\\
\dot{y}_i&=p_i^y\\
\dot{p}_i^x&=-\nu p_i^x-\frac{\partial V}{\partial x_i}=-\nu p_i^x-\sum_{j=1}^N\omega_{ij}(x_i-x_j)\\
\dot{p}_i^y&=-\nu p_i^y-\frac{\partial V}{\partial y_i}=-\nu p_i^y-\sum_{j=1}^N\omega_{ij}(y_i-y_j)
\end{aligned}\right.\;\; ,
\end{equation} $i=1,\ldots,N$, where $p_i^x$ and $p_i^y$ are the respective momenta
in the $x$ and the $y$ directions, assuming unit mass for each
particle. For this system, the configuration space is
$Q=\mathbb{R}^{2N}$ and $\mathbf{q}\in
Q$; an element in the cotangent bundle $\mathbf{z}\in T^*Q$ can be
written as
$\mathbf{z}=(\mathbf{q},\mathbf{p})=((\mathbf{x},\mathbf{y}),({\mathbf{p}^{x}},{\mathbf{p}^{y}}))$, where $\mathbf{p}^{x}=(p_1^x,\ldots, p_N^x)$ and $\mathbf{p}^{y}=(p_1^y,\ldots, p_N^y) $.

Our goal is to solve for and stabilize an arbitrary tensegrity structure
corresponding to a desired shape given by $[\mathbf{q}^e]$. The first step is to find the stresses $\omega_{ij}$ that make a given $[\mathbf{q}^e]$ an equilibrium shape of
the system. Recall that $\omega_{ij} = 0$ implies that nodes $i$ and $j$ are not connected.
$\mathbf{z}^e$ is an equilibrium of (\ref{eomfirstorder}) if and
only if $\mathbf{p}^{x}=\mathbf{p}^{y}=\mathbf{0}$ and
$\mathbf{q}^e = (\mathbf{x}^e, \mathbf{y}^e)$ is a critical point of the potential
$V(\mathbf{q})$.  Since $V(\mathbf{q})$ depends only on relative positions of nodes, then $\mathbf{q}^e$ is an equilibrium shape if and only if every $\mathbf{q} \in [\mathbf{q}^e]$ is an equilibrium shape. Let $ [\mathbf{z}^e]=( [\mathbf{q}^e],\mathbf{0},\mathbf{0})$, then $\mathbf{z}^e$ is an equilibrium of (\ref{eomfirstorder}) if an only if every $\mathbf{z}\in[\mathbf{z}^e]$ is an equilibrium of (\ref{eomfirstorder}).  

We now rewrite (\ref{potential1}) to find a more
exploitable relationship between the stresses $\omega_{ij}$ and
the critical points of $V(\mathbf{q})$. Using notations from
algebraic graph theory, we consider a tensegrity structure as an
undirected graph $G=(\mathcal{V},\mathcal{E})$, where
$\mathcal{V}$ is the set of nodes and $\mathcal{E}$ the set of
edges. Let $d_j$ be the degree of node $j$ (the number of edges at node $j$), then the Laplacian $L$
of the graph $G$ is the $N\times N$ matrix defined by
\begin{equation} \label{Laplacian}
L_{ij}=\left\{ \begin{array}{ll}
d_j & \textrm{if $i=j$} \\
-1 & \textrm{if $(i,j) \in \mathcal{E}$} \\
0 & \textrm{otherwise}.
\end{array} \right.
\end{equation}
In our setting, interconnection strengths are not identical from one edge to
the other but rather are weighted by the spring constants
$\omega_{ij}$. A tensegrity structure can then be viewed as an
undirected graph for which we define the weighted pseudo Laplacian
$\Omega$ by
\begin{displaymath} \label{stressmatrix}
\Omega_{ij}=\left\{ \begin{array}{ll}
\sum\limits_{j=1}^{N}\omega_{ij} & \textrm{if $i=j$}\\
-\omega_{ij} & \textrm{if $i\neq j$} .
\end{array} \right .
\end{displaymath}
This matrix $\Omega$ is called  the \emph{stress
matrix} by Connelly \cite{connelly1}. The stress matrix is an $N\times N$ symmetric matrix, and
the $N$-dimensional vector $\mathbf{1} = (1, \ldots, 1)^T$  is in the kernel of $\Omega$ (this last
property being true for all Laplacians). 

The potential (\ref{potential1}) can be rewritten as
\begin{equation}\label{potential2}
V(\mathbf{q})=\frac{1}{2}\mathbf{q}^{T}(\Omega\otimes
I_{2})\mathbf{q},
\end{equation}
where $I_{2}$ is the $2\times 2$ identity matrix and $\Omega\otimes
I_{2}$ is the $2n\times 2n$ block diagonal matrix $\begin{pmatrix}
\Omega & 0 \cr 0 & \Omega\end{pmatrix}$. This new form of the potential (\ref{potential2}) allows us to
derive a simple relationship between the choice of stresses
$\omega_{ij}$ and the equilibria of the model. The critical points
of (\ref{potential2}) and hence the equilibria $\mathbf{q}^e$ of
(\ref{eomfirstorder}) satisfy
\begin{equation}\label{equilibrium.equation}
\mathbf{q}^{T}(\Omega\otimes I_{2})=0.
\end{equation}
Since the stress matrix is symmetric, a placement $\mathbf{q}^e=(\mathbf{x}^e,\mathbf{y}^e)$
is a critical point of $V(\mathbf{q})$ if and only if $\mathbf{x}^e$
and $\mathbf{y}^e$ are each in the kernel of $\Omega$. Recall that
$\mathbf{1}$ is in the kernel of $\Omega$. Assuming that the nodes
are not all in a line, $\mathbf{x}^e,\mathbf{y}^e$ and
$\mathbf{1}$ are linearly independent. We can conclude that with
this model, a combination of cables and struts will have an
equilibrium if and only if $\mbox{rank}(\Omega)\leq N-3$. We
assume from now on that $N\geq 4$.

By choosing the stresses of the edges of the structure so that $\mbox{rank}(\Omega)=N-3$, the
kernel of the stress matrix $\Omega$ is exactly three dimensional.
We derive a method in Section~\ref{nonlinear} to compute $\Omega$ that fixes the desired shape of the equilibrium.  However, the computed $\Omega$ will not fix the scale of the equilibrium configuration for the model (\ref{eomfirstorder}). Indeed
if
$\mbox{ker}(\Omega)=\mbox{span}\{\mathbf{x}^e,\mathbf{y}^e,\mathbf{1}\}$
then $(\alpha \mathbf{x}^e ,\beta
\mathbf{y}^e,\mathbf{0},\mathbf{0})$ is also an equilibrium
$\forall \alpha,\beta \in\mathbb{R}$. For real tensegrities this
is not a problem because the cable and strut constraints preclude
the existence of any but the original equilibrium. In the virtual
setting, however, where the same constraints cannot be imposed, we
get a continuum of equilibria, which is not desirable. For
example, if the prescribed shape for the tensegrity is a square,
it will be the case that not only all squares but also all
rectangles will be equilibria.

In order to use tensegrity structures in our virtual setting, we
need to modify this model to isolate equilibria. In the next section we propose an augmented force model and show that an equilibrium shape (geometry and scale) is isolated.
We design the augmented model so that we can still make use of equation
(\ref{equilibrium.equation}), derived using the model
(\ref{force:model1}), that determines the geometry of the
tensegrity as a function of the parameters $\omega_{ij}$.

\subsection{Augmented model and smooth mapping from arbitrary shape to tensegrity}\label{nonlinear}

\qquad In the previous section we modeled edges as springs with zero rest length. Because in our virtual setting it is not possible to impose physical constraints on the edges, this model yields a continuum of equilibrium shapes for the system (\ref{eomfirstorder}), fixing only the geometry of the structure but not its scale. We now present an augmented model that allows us to isolate the shape geometry and scale.   We then derive a systematic and smooth method that prescribes tensegrity edge allocations and model parameters so that an arbitrary desired planar shape is an isolated equilibrium of the tensegrity dynamics.    In Section~\ref{sec:stability} we prove that the desired shape is an exponentially stable equilibrium of the tensegrity dynamics so that a vehicle network with control forces simulating forces internal to the derived tensegrity structure will converge exponentially to the desired shape.

In the augmented model, edges are modeled as springs with finite, nonzero rest length.
An edge is a cable when it is longer than its rest length, i.e., in tension; an edge is a strut when it is
shorter than its rest length, i.e., in compression. For two nodes $i, j$ we define
\begin{equation} \label{force:model2}
\vec f_{i\rightarrow j}=
\alpha_{ij}\omega_{ij}\frac{r_{ij}-l_{ij}}{r_{ij}}(\vec{q}_{i}-\vec{q}_{j})=-\vec
f_{j\rightarrow i}.
\end{equation}
Here $r_{ij}=\|\vec{q}_{i}-\vec{q}_{j}\|$ is the relative distance
between nodes $i$ and $j$, $l_{ij}>0$ is the rest length of the
spring that models the edge $ij$, $\omega_{ij}$ is the spring
constant from model (\ref{force:model1}), and $\alpha_{ij}$ is a
scalar parameter that fixes the spring constant of model
(\ref{force:model2}) for the edge $ij$.

We derive for this model an analogue of the stress
matrix $\Omega$. The potential energy of a tensegrity structure
with the forces induced by its edges modeled by
(\ref{force:model2}) is
\begin{equation} \label{potential4}
\tilde{V}(\mathbf{q})=\frac{1}{2}\sum_{i<j}\alpha_{ij}\omega_{ij}(r_{ij}-l_{ij})^{2}.
\end{equation}
With the same damping force used in the previous section, the equations
of motion are
\begin{equation}\label{eomfirstorder2}
\left\{ \begin{aligned}
\dot{x}_i&=p_i^x\\
\dot{y}_i&=p_i^y\\
\dot{p}_i^x&=-\nu p_i^x-\frac{\partial \tilde{V}}{\partial x_i}=-\nu p_i^x-\sum_{j=1}^{N}\tilde{\omega}_{ji}(x_i-x_j)\\
\dot{p}_i^y&=-\nu p_i^y-\frac{\partial \tilde{V}}{\partial y_i}=-\nu p_i^y-\sum_{j=1}^{N}\tilde{\omega}_{ji}(y_i-y_j)
\end{aligned}\right.\;\;
\end{equation}
$i=1,\ldots,N$, where $\widetilde{\omega}_{ij}$ is given by
\begin{equation} \label{stressnonlinear}
\widetilde{\omega}_{ij}(\mathbf{x},\mathbf{y})=\alpha_{ij}\omega_{ij}(1-\frac{l_{ij}}{r_{ij}}).
\end{equation}

As before, our goal is to solve for and stabilize a tensegrity structure corresponding to a desired shape $[\mathbf{q}^e]$.  First we compute the $\tilde{\omega}_{ij}$ so that $[\mathbf{q}^e]$ is an equilibrium shape of the
system (\ref{eomfirstorder2}). This requires finding the relationship between the choice of the
$\omega_{ij},\alpha_{ij}$ and $l_{ij}$'s and the equilibria of the
system (\ref{eomfirstorder2}).    Since $\tilde V(\mathbf{q})$ depends only on relative positions of nodes, then $\mathbf{q}^e$ is an equilibrium shape if and only if every $\mathbf{q} \in [\mathbf{q}^e]$ is an equilibrium shape.  As in the previous model,
$\mathbf{z}^e$ is an equilibrium of (\ref{eomfirstorder2}) if and
only if $\mathbf{p}^{x}=\mathbf{p}^{y}=\mathbf{0}$ and
$\mathbf{q}^e = (\mathbf{x}^e, \mathbf{y}^e)$ is a critical point of the potential
$\tilde V(\mathbf{q})$. We likewise define $[\mathbf{z}^e]=( [\mathbf{q}^e],\mathbf{0},\mathbf{0})$, then $\mathbf{z}^e$ is an equilibrium of (\ref{eomfirstorder2}) if an only if every $\mathbf{z}\in[\mathbf{z}^e]$ is an equilibrium of (\ref{eomfirstorder2}).  
 The critical points of (\ref{potential4}) are
given by
\begin{equation}\label{criticalpoints}
\begin{aligned}
\sum_{j=1}^{N}\alpha_{ij}\omega_{ij}(\vec{q}_j-\vec{q}_i)\left(1-\frac{l_{ij}}{r_{ij}}\right)=0,
&{}\;\; i=1,\ldots,N.
\end{aligned}
\end{equation}
From (\ref{criticalpoints}), an analogue of the stress matrix
$\Omega$ is constructed. The new stress matrix is not a constant
matrix and depends on the relative distances  $r_{ij}$ between pairs of
nodes.  The entries of the new stress
matrix are given by
\begin{displaymath}
\widetilde\Omega_{ij}(\mathbf{x},\mathbf{y})=\left\{
\begin{array}{ll}
\sum\limits_{j=1}^{N}\widetilde{\omega}_{ij}(\mathbf{x},\mathbf{y})
& \textrm{if $i=j$} \\
-\widetilde{\omega}_{ij}(\mathbf{x},\mathbf{y}) & \textrm{if
$i\neq j$,}
\end{array} \right.
\end{displaymath}
where the $\widetilde{\omega}_{ij}$ is the new stress for the edge $ij$ given by (\ref{stressnonlinear}).
The vector $\mathbf{1}$ is also in the kernel of
$\widetilde{\Omega}$, $\forall$ $(\mathbf{x},\mathbf{y})\in
\mathbb{R}^N\times\mathbb{R}^N$.

This new stress matrix can now be used to characterize the
equilibria of (\ref{eomfirstorder2}). To make the placement
$\mathbf{q}^e=(\mathbf{x}^e,\mathbf{y}^e)$ a critical point of $\tilde V(\mathbf{q})$, we pick (if possible) the parameters
$\alpha_{ij},\omega_{ij},l_{ij}$ so that
\begin{equation}\label{eqcond}
\begin{aligned}
\widetilde{\Omega}(\mathbf{x}^e,\mathbf{y}^e)\mathbf{x}^e & =0\\
\widetilde{\Omega}(\mathbf{x}^e,\mathbf{y}^e)\mathbf{y}^e & =0.
\end{aligned}
\end{equation}
As a first step we choose the parameters $\alpha_{ij}$ and
$l_{ij}$ for all $i,j$ so that
$\widetilde{\Omega}(\mathbf{x}^e,\mathbf{y}^e)=\Omega$. In order to make $\widetilde{\omega}_{ij}(\mathbf{x}^e,\mathbf{y}^e)=\omega_{ij}$, we choose $\alpha_{ij}, l_{ij}$ such that $\alpha_{ij}(1-\frac{l_{ij}}{r_{ij}^{e}})=1$, where $r_{ij}^{e}$ is the value of $r_{ij}$ at $[\mathbf{q}^e]$. This last equation is solved by picking
\begin{equation}\label{alphal}
\begin{aligned}
\alpha_{ij}&=\frac{\pi}{\arctan\omega_{ij}}\\
l_{ij}&=r_{ij}^e\left(1-\frac{1}{\pi}\arctan\omega_{ij}\right).
\end{aligned}
\end{equation}
If edge $ij$ is a strut, then $\omega_{ij}<0$, and we have from equation (\ref{alphal}) that $\alpha_{ij}<0$ which makes the spring constant $\alpha_{ij}\omega_{ij}>0$ and $l_{ij}>r_{ij}^e$. If edge $ij$ is a cable then $\omega_{ij}>0$, and we have from equation (\ref{alphal}) that $\alpha_{ij}>0$ which makes the spring constant $\alpha_{ij}\omega_{ij}>0$ and $l_{ij}<r_{ij}^e$. 
The choice of $l_{ij}$ and $\alpha_{ij}$ is not unique, but this choice makes the vector field (\ref{eomfirstorder2})  a $C^{\infty}$ map of $\omega_{ij}$ and $r_{ij}^e$. This result will be critical in Section \ref{shape} to prove stability for the time-varying control law that changes the shape of the formation from any initial shape to any final desired shape.

We now show that the parameters $\omega_{ij}$ can be found
independently of parameters $\alpha_{ij}$ and $l_{ij}$, so that
$\mbox{ker}(\Omega)=\mbox{span}\{\mathbf{x}^e,\mathbf{y}^e,\mathbf{1}\}$,
and such that the nonzero eigenvalues of $\Omega$ are all
positive. This makes the equilibrium $\mathbf{z}^e=(\mathbf{x}^e,
\mathbf{y}^e,\mathbf{0},\mathbf{0})$ an isolated minimum of the
potential (modulo rigid transformations), i.e., our choices ensure
that we have the right combination of struts and cables to make
$\mathbf{q}^e=(\mathbf{x}^e, \mathbf{y}^e)$ (and therefore any $\mathbf{q} \in [\mathbf{q}^e]$) a tensegrity
structure. $\Omega$ is a symmetric matrix, hence it has only real
eigenvalues and it can be diagonalized using an orthonormal basis.
As mentioned previously, if the case when all the nodes are in a
line is excluded, then $\mathbf{x}^e,\mathbf{y}^e$ and
$\mathbf{1}$ are linearly independent. We can complete these three
vectors with $N-3$ others and obtain a basis of $\mathbb{R}^N$.
Then applying the Gram-Schmidt procedure to those vectors yields
an orthonormal basis $(\mathbf{v}_1,\ldots,\mathbf{v}_N)$ for
$\mathbb{R}^N$ that satisfies
\begin{equation*}
\mbox{span}\{\mathbf{v}_1,\mathbf{v}_2,\mathbf{v}_3\}=\mbox{span}\{\mathbf{x}^e,\mathbf{y}^e,\mathbf{1}\}.
\end{equation*}
Now we define the $N\times N$ diagonal matrix $D=\mbox{diag}(0,0,0, d_4, \ldots, d_N)$, where $d_i>0,$ $\forall i\in \{4, \cdots, N\}$, and the orthonormal $N\times N$
matrix $\Lambda=(\mathbf{v}_1  \cdots 
\mathbf{v}_N)$. If we compute $\Lambda D\Lambda^{T}$
we have a symmetric positive semi-definite matrix with its kernel
equal to $\mbox{span}\{\mathbf{x}^e,\mathbf{y}^e,\mathbf{1}\}$.
Setting $\Omega=\Lambda D\Lambda^{T}$ determines the values of
stresses $\omega_{ij}$ that make the desired shape $[\mathbf{q}^e]$
 a tensegrity structure.    We prove stability of this equilibrium shape in Section~\ref{sec:stability}.

The choice of eigenvalues $D$ and eigenvectors $\Lambda$ for the stress matrix is not unique.
We investigate next, through an example, the influence of the choice of the eigenvalues and eigenvectors of the stress matrix on the resulting interconnection topology that achieves the given desired shape.

\subsection{Example}

We illustrate here the computation of the stress matrix using the method developed in Section \ref{nonlinear} and examine resulting interconnection topologies. It is known that stable tensegrity structures with $N$ nodes require at least $2N-2$ edges \cite{connelly2}. Connelly has proven a result that provides a means to systematically design stable planar tensegrities in the shape of any strictly convex polygon \cite{connelly1}. His method is designed to yield tensegrity structures with minimal number of edges. Our method allows us to generate  tensegrity structures of any shape, convex or non-convex, but often yields interconnection topologies with number of edges greater than the proven lower bound. However, as we illustrate , it is possible in our method to use the freedom of choice in the eigenvalues and eigenvectors of the stress matrix to reduce the number of edges. 

\begin{figure}[h]
\begin{center}
\includegraphics[width=2in]{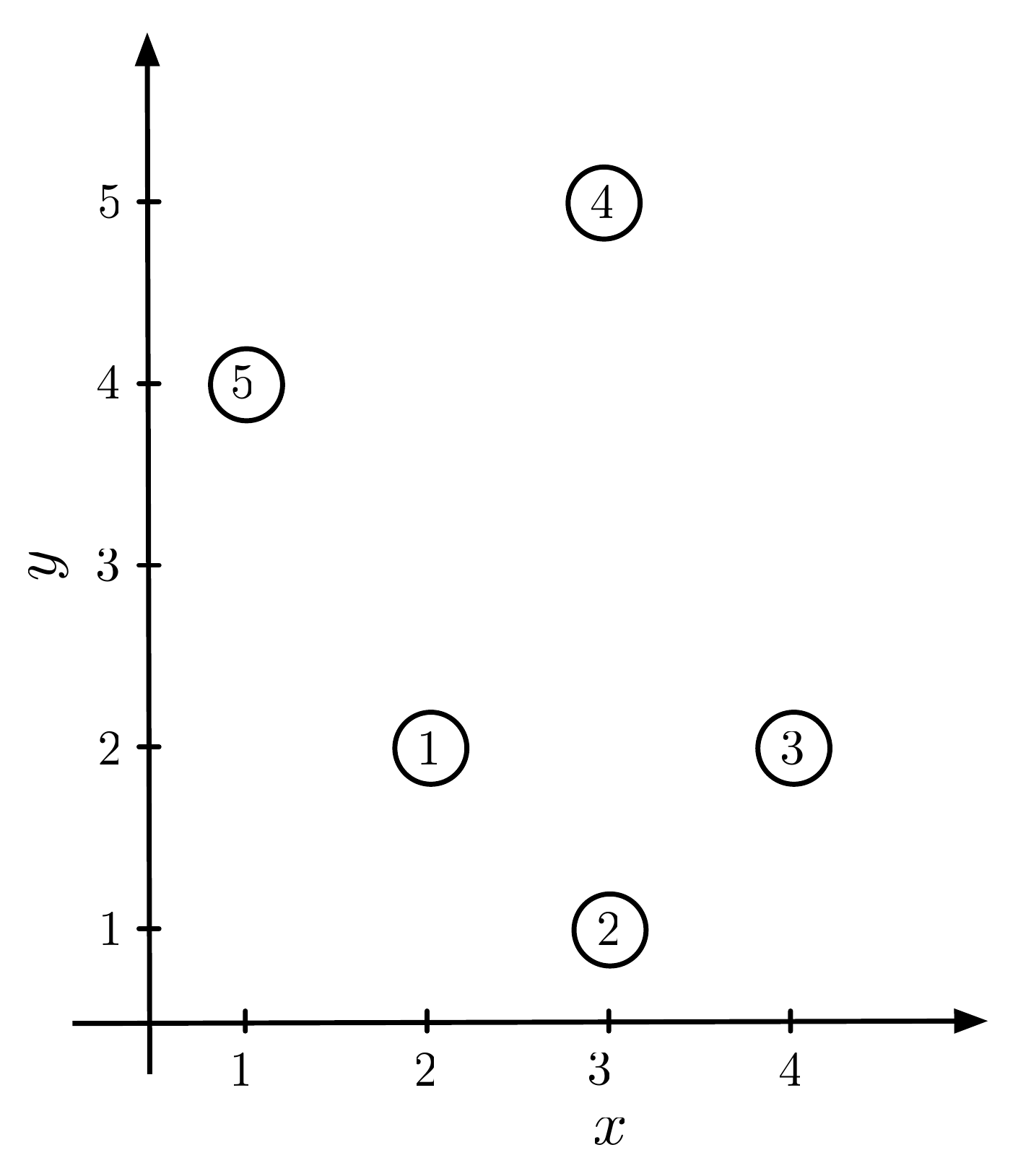}
\end{center}
\caption{Five node configuration given by $\mathbf{x}^e=( 2 , 3 , 4 , 3 ,1),
\mathbf{y}^e=( 2 , 1 , 2 , 5 , 4)$. The labels on the nodes define the ordering.} \label{fivenodesex}
\end{figure}

\begin{figure}[h]
\begin{center}
(a)\includegraphics[width=2in]{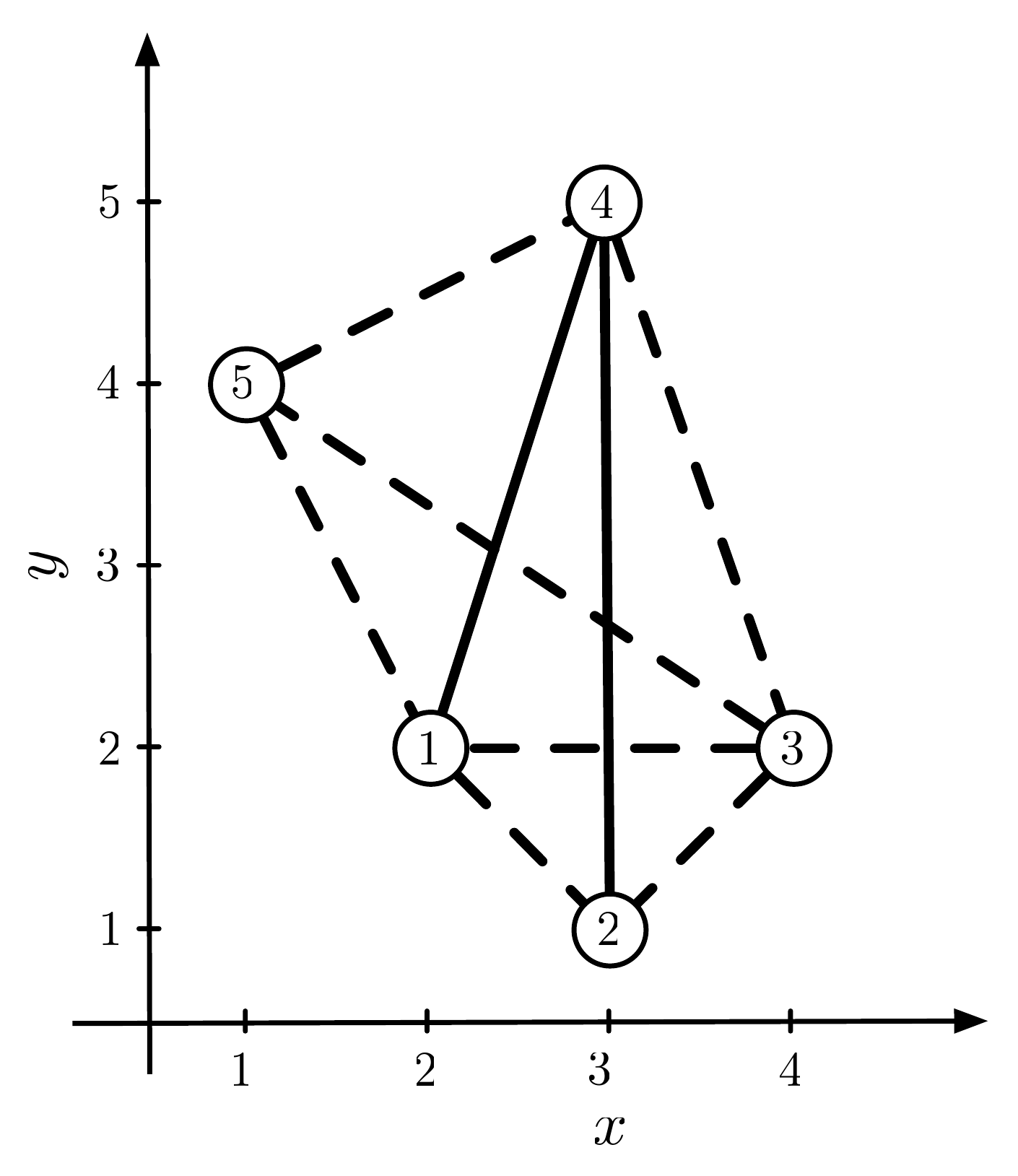}
(b)\includegraphics[width=2in]{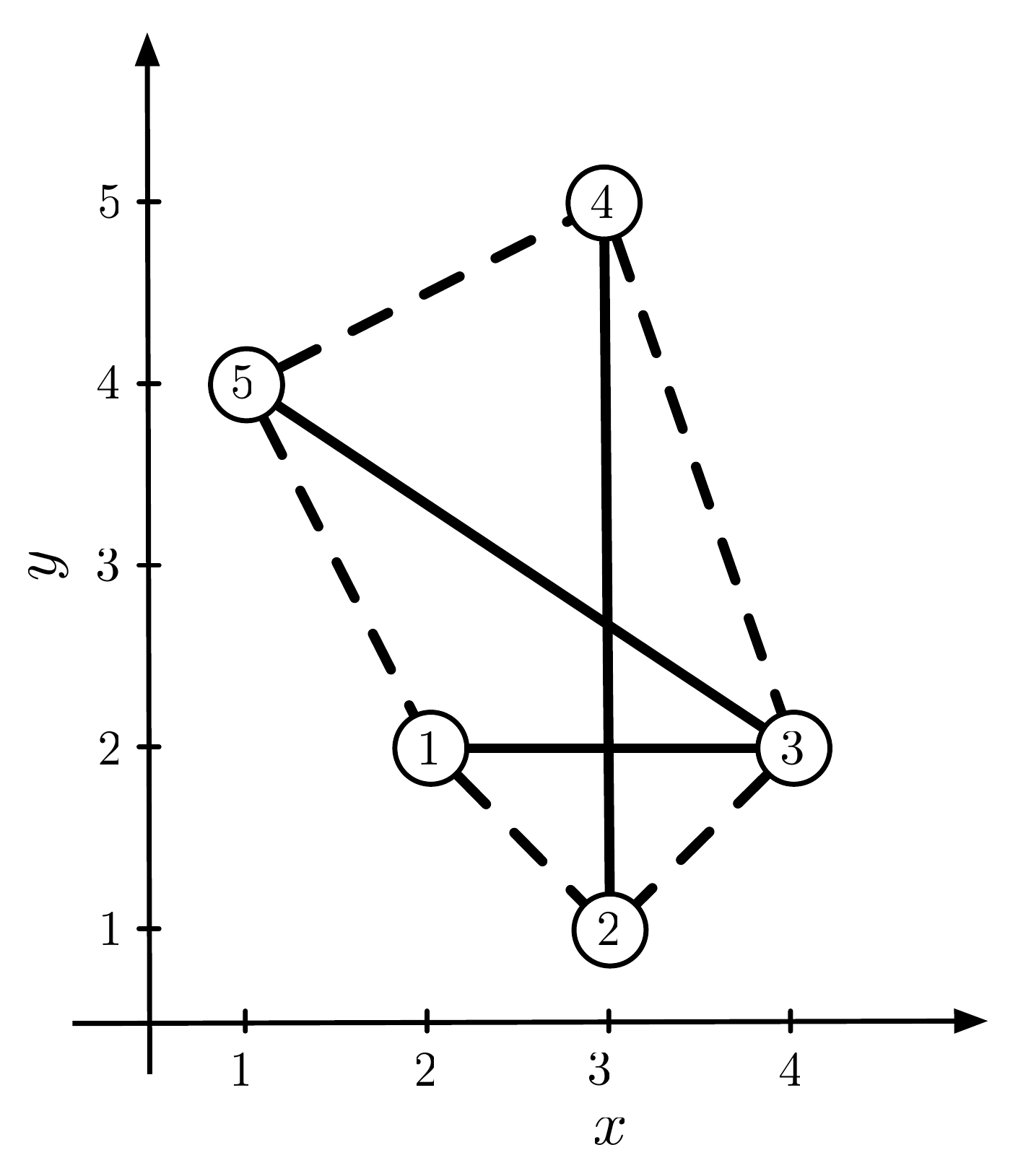}
\end{center}
\caption{Tensegrity structures generated from the stress matrices (a) $\Omega_1$, (b) $\Omega_2$.  Solid lines are struts and dashed lines are cables} \label{fivenodestens}
\end{figure}

Consider the network configuration given by 
\begin{equation}\label{irregpent}
\mathbf{q}^e=\begin{pmatrix}\mathbf{x}^e\cr\mathbf{y}^e\end{pmatrix}=\begin{pmatrix} 2 & 3 & 4 & 3 &1\cr 2 & 1 & 2 & 5 & 4\end{pmatrix}
\end{equation}and plotted in Figure \ref{fivenodesex}. 
Following the method from Section \ref{nonlinear} we take $D_1=\mbox{diag}(0,0,0,1,1)$ and
\begin{equation*}
\Lambda_1=\begin{pmatrix}\frac{\sqrt{5}}{5} & -\frac{3}{\sqrt{130}} & -\frac{\sqrt{182}}{39} & -\frac{\sqrt{518}}{222} & \frac{10\sqrt{74}}{111} \cr \frac{\sqrt{5}}{5} & \frac{2}{\sqrt{130}} & -\frac{\sqrt{182}}{26} & -\frac{\sqrt{518}}{37} & -\frac{3\sqrt{74}}{74} \cr \frac{\sqrt{5}}{5} & \frac{7}{\sqrt{130}} & -\frac{\sqrt{182}}{273} & \frac{43\sqrt{518}}{1554} & -\frac{2\sqrt{74}}{111} \cr \frac{\sqrt{5}}{5} & \frac{2}{\sqrt{130}} & -\frac{31\sqrt{182}}{546} & -\frac{13\sqrt{518}}{777} & \frac{5\sqrt{74}}{222} \cr \frac{\sqrt{5}}{5} & \frac{-8}{\sqrt{130}} & \frac{2\sqrt{182}}{182} & \frac{3\sqrt{518}}{259} & -\frac{2\sqrt{74}}{37}\end{pmatrix},
\end{equation*}
where the columns of $\Lambda_1$ constitute a basis of orthonormal eigenvectors obtained with the Gram-Schmidt procedure.   The stress matrix $\Omega_1=\Lambda_1 D_1\Lambda_1^T$ is given by
\begin{equation}\label{stress1}
\Omega_1=\begin{pmatrix} \frac{11}{18} & -\frac{1}{3} & -\frac{1}{18} & \frac{1}{9} & -\frac{1}{3} \cr -\frac{1}{3} & \frac{1}{2} & -\frac{1}{3} & \frac{1}{6} & 0 \cr -\frac{1}{18} & -\frac{1}{3} & \frac{53}{126} & -\frac{17}{63} & -\frac{5}{21} \cr \frac{1}{9} & \frac{1}{6} & -\frac{17}{63} & \frac{23}{126} & -\frac{4}{21}   \cr -\frac{1}{3} & 0 & \frac{5}{21} & -\frac{4}{21} & \frac{2}{7} \end{pmatrix}.
\end{equation}
The tensegrity structure corresponding to this stress matrix, plotted in Figure \ref{fivenodestens}(a) has an  interconnection topology requiring 9 edges. To reduce the number of edges, we manipulate our choice of $D$ and $\Lambda$ so that $\Omega = \Lambda D \Lambda^T$ has entries identically equal to zero.  Consider, for example, setting $D_2=\mbox{diag}\begin{pmatrix} 0 & 0 & 0 & \frac{60}{253} & \frac{30}{253}\end{pmatrix}$ and
\begin{equation*}
\Lambda_2=\begin{pmatrix} 0 & 0 & \frac{\sqrt{14}}{6} & \frac{\sqrt{2}}{3} & -\frac{\sqrt{14}}{6}\cr \frac{1}{\sqrt{42}} & -\frac{2}{\sqrt{21}} & \frac{\sqrt{14}}{7} & -\frac{\sqrt{2}}{2} & 0 \cr \frac{4}{\sqrt{42}} & -\frac{2}{\sqrt{21}} & \frac{\sqrt{14}}{42} & \frac{\sqrt{2}}{3} & \frac{5\sqrt{14}}{42} \cr \frac{5}{\sqrt{42}} & \frac{2}{\sqrt{21}} & -\frac{\sqrt{14}}{21} & -\frac{\sqrt{2}}{6} & -\frac{2\sqrt{14}}{21} \cr 0 & \frac{3}{\sqrt{21}} & \frac{\sqrt{14}}{7} & 0 & \frac{\sqrt{14}}{7} \end{pmatrix}.
\end{equation*}
We compute the stress matrix:
\begin{equation}\label{stress2}
\Omega_2=\begin{pmatrix} \frac{25}{253} & -\frac{20}{253} & \frac{5}{253} & 0 & -\frac{10}{253} \cr -\frac{20}{253} & \frac{30}{253} & -\frac{20}{253} & \frac{10}{253} & 0 \cr \frac{5}{253} & -\frac{20}{253} & \frac{135}{1771} & -\frac{80}{1771} & \frac{50}{1771} \cr 0 & \frac{10}{253} & -\frac{80}{1771} & \frac{50}{1771} & -\frac{40}{1771}   \cr -\frac{10}{253} & 0 & \frac{50}{1771} & -\frac{40}{1771} & \frac{60}{1771} \end{pmatrix}.
\end{equation}
The tensegrity corresponding to this stress matrix, plotted in Figure \ref{fivenodestens}(b) has an interconnection topology  requiring only $8=2\times 5-2$ edges. The number of edges cannot be reduced further since we have reached the lower bound proven in \cite{connelly2}. 

As the number of nodes increases, it becomes harder to systematically find a combination of eigenvalues and eigenvectors that yields tensegrities with minimal number of edges.  
However, as the above example illustrates, it may be possible to manipulate the choices of $D$ and $\Lambda$ to reduce the number of connections in a tensegrity derived using our method.  

Our approach has significant advantages because it is systematic and smooth.   We exploit this in Section \ref{shape}, where we define a smooth parameterization of the control law (\ref{eomfirstorder2}), creating a systematic framework for smooth reconfiguration of tensegrity structures between arbitrary planar shapes. The choice of parameters $\alpha_{ij}, l_{ij}$ and $\omega_{ij}$ as defined in Section \ref{nonlinear} is crucial in defining this smooth parameterization. 

\section{Stability Analysis}\label{sec:stability}

In Section \ref{nonlinear} we presented an augmented model for the
forces along the edges of a tensegrity and a systematic method to
determine the parameters of the model (\ref{force:model2}) that
make any desired shape in the plane an equilibrium of
(\ref{eomfirstorder2}). We now prove for the parameters
$\alpha_{ij},l_{ij},\omega_{ij}$ satisfying (\ref{eqcond}), that the
equilibrium set $[\mathbf{z}^e]=([\mathbf{q}^e],\mathbf{0},\mathbf{0})$ of (\ref{eomfirstorder2})
is  locally exponentially stable. We first
show that  $[\mathbf{q}^e]$ is
an isolated critical point of the potential
$\tilde V(\mathbf{q})$ and $[\mathbf{z}^e]$ an isolated equilibrium set of
(\ref{eomfirstorder2}). We then
prove local exponential stability for
$[\mathbf{z}^e]$ using the linearization of (\ref{eomfirstorder2}).

With the parameters $\alpha_{ij}, l_{ij}$ and $\omega_{ij}$ satisfying
(\ref{eqcond}), $\mathbf{q}^e=(\mathbf{x}^e, \mathbf{y}^e)$ is a critical point of the
potential $\tilde V(\mathbf{q})$. The potential $\tilde V(\mathbf{q})$ only
depends on the relative distances between the nodes. This implies
that $\tilde V(\mathbf{q})$ and hence the Lagrangian of the system are
invariant under the action of the Lie group $SE(2)$ on the
configuration space $Q$. Given these symmetries, critical points
of $\tilde V(\mathbf{q})$ and hence equilibria of (\ref{eomfirstorder2})
are only isolated modulo $SE(2)$ transformations. We now show that
$[\mathbf{q}^e]$ is an isolated set of critical points of $\tilde V(\mathbf{q})$ by computing
$\delta^2 \tilde V(\mathbf{q})$ and proving that when the second variation
of $\tilde V(\mathbf{q})$ is evaluated at $\mathbf{q}^e$, we have a
positive definite matrix except in the symmetry directions
$SE(2)$. $\delta^2 \tilde V(\mathbf{q}^e)$ is given by
\begin{equation}\label{secondvariation}
\delta^2 \tilde V(\mathbf{q}^e)=\begin{pmatrix} \Omega+L_{\omega
x}(\mathbf{q}^e) & L_{\omega xy}(\mathbf{q}^e)\cr L_{\omega
xy}(\mathbf{q}^e) & \Omega+L_{\omega y}(\mathbf{q}^e)\end{pmatrix}
\end{equation} where $\Omega$ is the stress matrix derived in Section 2 and the $ij$th element of each block matrix is
\begin{displaymath}
L_{\omega x}(i,j)=\left\{\begin{array}{ll}
-\alpha_{ij}\omega_{ij}\frac{(x_{i}-x_{j})^{2}l_{ij}}{r_{ij}^{3}}
& \textrm{if $i\neq j$}\\
\sum\limits_{j=1, j\neq
i}^{N}\alpha_{ij}\omega_{ij}\frac{(x_{i}-x_{j})^{2}l_{ij}}{r_{ij}^{3}}
& \textrm{if $i=j$} \end{array} \right.
\end{displaymath}
\begin{displaymath}
L_{\omega y}(i,j)=\left\{\begin{array}{ll}
-\alpha_{ij}\omega_{ij}\frac{(y_{i}-y_{j})^{2}l_{ij}}{r_{ij}^{3}}
& \textrm{if $i\neq j$}\\
\sum\limits_{j=1, j\neq
i}^{N}\alpha_{ij}\omega_{ij}\frac{(y_{i}-y_{j})^{2}l_{ij}}{r_{ij}^{3}}
& \textrm{if $i=j$} \end{array} \right.
\end{displaymath}
and
\begin{displaymath}
L_{\omega xy}(i,j)=\left\{\begin{array}{ll}
-\alpha_{ij}\omega_{ij}\frac{(x_{i}-x_{j})(y_{i}-y_{j})l_{ij}}{r_{ij}^{3}}
& \textrm{if $i\neq j$}\\
\sum\limits_{j=1, j\neq
i}^{N}\alpha_{ij}\omega_{ij}\frac{(x_{i}-x_{j})(y_{i}-y_{j})l_{ij}}{r_{ij}^{3}}
& \textrm{if $i=j$.} \end{array} \right.
\end{displaymath}

\begin{positive}\label{pos}$\delta^2 \tilde V(\mathbf{q}^e)$ is a positive
semi-definite matrix.
\end{positive}
\textbf{Proof:} See proof in Appendix \ref{lemma1}.

\begin{kernel}\label{ker} The kernel of $\delta^2 \tilde V(\mathbf{q}^e)$ is equal to
\[
\mbox{span} \left\{ \left( \begin{array}{c} \mathbf{1} \\
\mathbf{0}\end{array} \right),  \left( \begin{array}{c} \mathbf{0} \\
\mathbf{1}\end{array} \right),
\left( \begin{array}{c} -\mathbf{y}^e \\
\mathbf{x}^e\end{array} \right) \right\}.
\]
\end{kernel}
\textbf{Proof:} See proof in Appendix \ref{lemma2}.

\begin{isolated}\label{isolated} $[\mathbf{z}^e]=([\mathbf{q}^e], \mathbf{0},
\mathbf{0})$  is an isolated equilibrium set of (\ref{eomfirstorder2}).
\end{isolated}

\textbf{Proof:} By Lemma \ref{ker}, three linearly independent eigenvectors for the three zero eigenvalues are $\left\{ \left( \begin{array}{c} \mathbf{1} \\
\mathbf{0}\end{array} \right),  \left( \begin{array}{c} \mathbf{0} \\
\mathbf{1}\end{array} \right),
\left( \begin{array}{c} -\mathbf{y}^e \\
\mathbf{x}^e\end{array} \right) \right\}$. These vectors
correspond respectively to symmetries of translation along the
$x$-axis and the $y$-axis and of rotation about the origin. In
addition, Lemma \ref{pos} guarantees that all other eigenvalues of
$\delta^2V(\mathbf{q}^e)$ are strictly positive. Combining Lemma
\ref{pos} and Lemma \ref{ker} concludes the proof of the theorem.
$\quad \square$

We note that the equilibrium $[\mathbf{z}^e]$ is not necessarily a
unique equilibrium.

\begin{stable}\label{stable}
$[\mathbf{z}^e]=([\mathbf{q}^e],\mathbf{0},\mathbf{0})$ is a locally exponentially stable equilibrium set of
(\ref{eomfirstorder2}).
\end{stable}

\textbf{Proof:} We prove local exponential stability for the
isolated equilibrium set $[\mathbf{z}^e]$ using linearization of
(\ref{eomfirstorder2}). The Jacobian of
(\ref{eomfirstorder2}) evaluated at
$\mathbf{z}^e=(\mathbf{x}^e,\mathbf{y}^e,\mathbf{0},\mathbf{0})$
is
\begin{equation}\label{jacobian}
\mbox{D}g(\mathbf{z}^e)=\begin{pmatrix} 0_{2n} & I_{2n}\cr
-\delta^2 \tilde V(\mathbf{q}^e) & -\nu I_{2n}\end{pmatrix},
\end{equation}
where the vector field $g$ represents the right hand side of (\ref{eomfirstorder2}).  We show that this matrix is negative semi-definite and the three
zero eigenvalues correspond to the symmetry directions $SE(2)$.
The linearization $\mbox{D}g(\mathbf{z}^e)$ can be written as
\begin{equation*}
\mbox{D}g(\mathbf{z}^e)=\begin{pmatrix}
\frac{1}{\nu}\delta^2 \tilde V(\mathbf{q}^e) & I_{2n}\cr
-\delta^2 \tilde V(\mathbf{q}^e) & -\nu
I_{2n}\end{pmatrix}+\begin{pmatrix}
-\frac{1}{\nu}\delta^2 \tilde V(\mathbf{q}^e) & 0_{2n}\cr 0_{2n} &
0_{2n}\end{pmatrix}=: B_1+B_2.
\end{equation*}
By Lemma \ref{pos}, $B_2$ is negative semi-definite. To show that
$B_1$ is also negative semi-definite, we proceed to two changes of
basis given by the following two invertible matrices:
\begin{align*}
P_1&= \begin{pmatrix} \frac{1}{\nu}I_{2n} & \frac{1}{\nu}I_{2n}\cr
0_{2n} & I_{2n}\end{pmatrix}\\
P_2&= \begin{pmatrix} 0_{2n} & I_{2n}\cr I_{2n} &
\frac{1}{\nu}\delta^2 \tilde V(\mathbf{q}^e)\left(\frac{1}{\nu}\delta^2 \tilde V(\mathbf{q}^e)-\nu
I_{2n}\right)^{-1}\end{pmatrix}.
\end{align*}
The representation of $B_1$ in the new basis is 
\begin{equation*}
P_2^{-1}P_1^{-1}D_1 P_1
P_2=\begin{pmatrix}\frac{1}{\nu}\delta^2 \tilde V(\mathbf{q}^e)-\nu I_{2n}
&0_{2n} \cr 0_{2n} & 0_{2n}\end{pmatrix}.
\end{equation*}
This matrix is negative semi-definite for all $
\nu>\sqrt{\lambda_{max}(\delta^2 \tilde V(\mathbf{q}^e))}$, where
$\lambda_{max}\left(\delta^2 \tilde V(\mathbf{q}^e)\right)$ is the largest
eigenvalue of $\delta^2 \tilde V(\mathbf{q}^e)$. We now show by direct
computation that the zero eigenvalues correspond to the $SE(2)$
symmetries:
\begin{equation*}
\begin{pmatrix} \mathbf{0}_{2n} & I_{2n}
\cr -\delta^2 \tilde V(\mathbf{q}^e) & -\nu
I_{2n}\end{pmatrix}\begin{pmatrix}\mathbf{x} \cr \mathbf{y} \cr
\mathbf{p}_x\cr
\mathbf{p}_y\end{pmatrix}=\begin{pmatrix}\mathbf{0} \cr \mathbf{0}
\cr \mathbf{0}\cr \mathbf{0}\end{pmatrix}
\end{equation*}
if and only if
\begin{equation*}
\begin{pmatrix}\mathbf{p}_x \cr \mathbf{p}_y \cr
-\delta^2  \tilde V(\mathbf{q}^e)\begin{pmatrix}\mathbf{x}\cr\mathbf{y}\end{pmatrix}-\nu
I_{2n}\begin{pmatrix}\mathbf{p}_x\cr \mathbf{p}_y\end{pmatrix}
\end{pmatrix}=\begin{pmatrix}\mathbf{0} \cr \mathbf{0} \cr
\mathbf{0}\cr \mathbf{0}\end{pmatrix}
\end{equation*}
if and only if
\begin{equation*}
\begin{pmatrix}\mathbf{x} \cr \mathbf{y} \cr
\mathbf{p}_x\cr
\mathbf{p}_y\end{pmatrix}\in\mbox{span} \left\{ \begin{pmatrix}\mathbf{1}
\cr \mathbf{0} \cr \mathbf{0}\cr
\mathbf{0}\end{pmatrix},\begin{pmatrix}\mathbf{0} \cr \mathbf{1}
\cr \mathbf{0}\cr
\mathbf{0}\end{pmatrix},\begin{pmatrix}-\mathbf{y}^e \cr
\mathbf{x}^e \cr \mathbf{0}\cr \mathbf{0}\end{pmatrix}.\right\}\;\;\;
\square
\end{equation*}

Local asymptotical stability can also be proved using the total
energy of the system as a Lyapunov function. In the next section
we utilize the exponential stability result
to present a time-dependent control law, defined as a parameterization of the control law used in (\ref{eomfirstorder2}), that enables a tensegrity to reconfigure itself between arbitrary desired planar shapes.

\section{Change of Shape}\label{shape}

The results from the previous sections give a framework that allow us to
stabilize any planar formation by mapping each vehicle to a node of
the tensegrity structure and controlling them with the forces
induced by the tensegrity's edges modeled by (\ref{force:model2}). We now take this result a
step further and present a control law for the nodes that enables a well-behaved reconfiguration between arbitrary shapes. The control law is designed so that the nodes follow a smooth path in the space of tensegrities.  Given a desired starting and an ending shape $[\mathbf{q}_0^e]$ and $[\mathbf{q}_f^e]$, we design a smooth path in shape space $[\mathbf{q}^e](t)$ such that $[\mathbf{q}^e](0)=[\mathbf{q}_0^e]$ and $[\mathbf{q}^e](\tau)=[\mathbf{q}_f^e]$, and we develop a method to deduce the necessary variations of parameters $\alpha_{ij}(t),l_{ij}(t)$ and $\omega_{ij}(t)$ to make a tensegrity follow that path in shape space.  The smoothly time-varying parameters define a smoothly time-varying control law.  We study the resulting controlled time-varying dynamical system and use results from Lawrence and Rugh \cite{Lawrence} on nonlinear systems with slowly varying inputs to prove the system is well behaved. By well behaved we mean boundedness during the trajectory, i.e. $\|[\mathbf{q}](t)-[\mathbf{q}^e](t)\|$ is bounded for $t \in [0,\tau]$, as well as convergence to the final shape, i.e., $[\mathbf{q}](t)\rightarrow [\mathbf{q}^e_f]$ as $t\rightarrow\infty$.

\subsection{Parameterized Control Law}

 Designing a path requires us to define pairings between the nodes of the initial configuration and the nodes of the final configuration. With an eye toward preventing collisions between vehicles and minimizing energy input to vehicles, we make a choice so that (1) the planned trajectories of  any two distinct nodes do not intersect and (2) the planned total distance travelled by all the nodes is minimized. We propose a solution where each vehicle travels on a straight line. Because of the symmetry in the system, the configuration can drift (rotate or translate) and therefore we cannot guarantee collision avoidance and minimal energy consumption. However simulations in section \ref{simulations} show good performence of the system with the well designed shape tracking laws we are now presenting.

\begin{figure}[h]
\begin{center}
\includegraphics[width=2in]{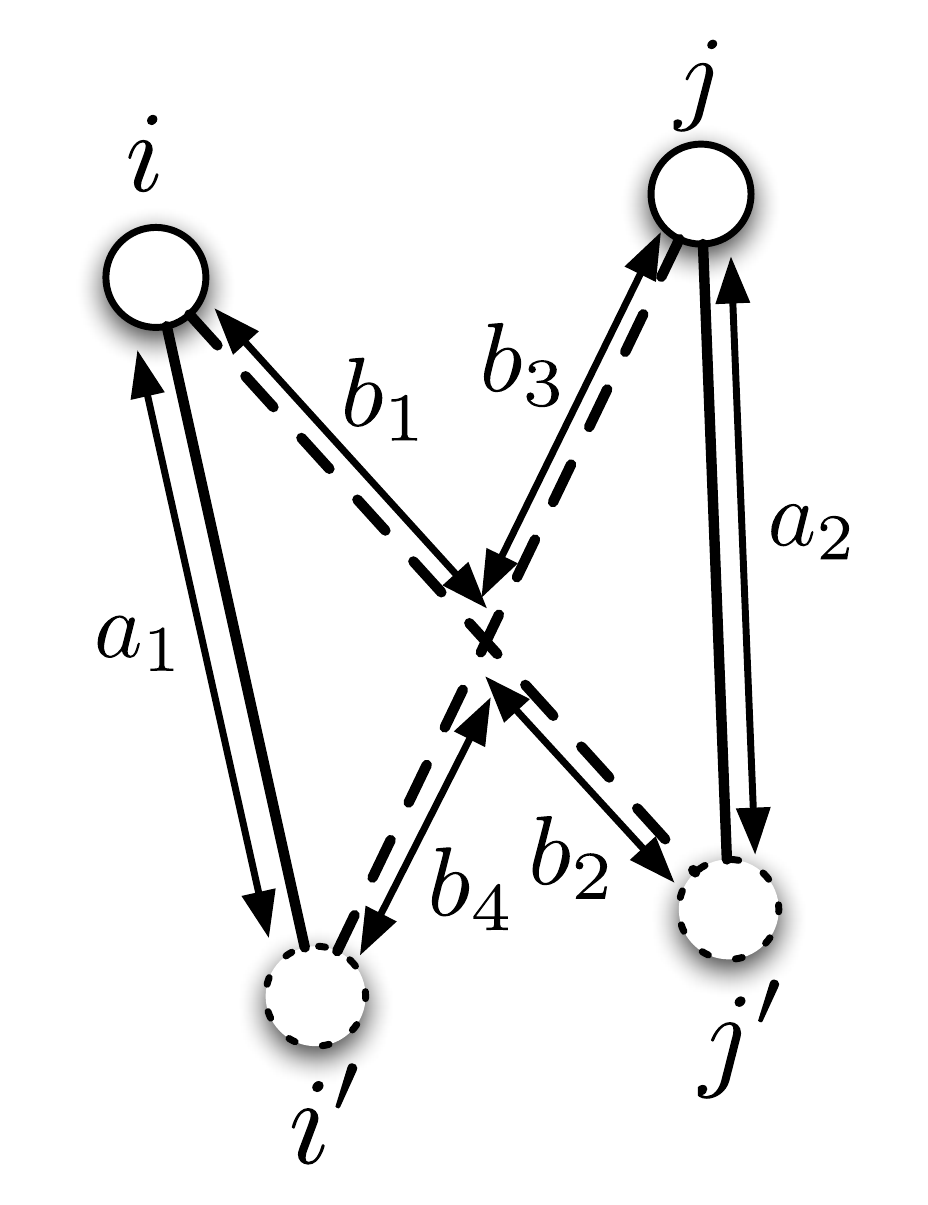}
\end{center}
\caption{Two possible pairings between two starting nodes $i,j$
and two ending nodes $i',j'$. For the intersecting trajectories
(dotted lines) the total distance travelled by the nodes is equal to
$b_1+b_2+b_3+b_4$. For the non-intersecting trajectories (solid
lines) the total distance travelled by the nodes is equal to
$a_1+a_2<b_1+b_2+b_3+b_4$.} \label{trajcross}
\end{figure}

\paragraph{Path planning:} To plan a path, we select representative configurations $\mathbf{q}_0\in[\mathbf{q}_0^e]$ and $\mathbf{q}_f\in[\mathbf{q}_f^e]$ that have the same centroid. We choose the orientations and the pairings between initial and final nodes to minimize
\begin{equation*}
\sum_{i=1}^N\|\vec{q}_{0i}-\vec{q}_{fi}\|.
\end{equation*}
This choice minimizes the planned total distance travelled by all the vehicles and ensures that no planned trajectories intersect.
\begin{pairing}\label{pairing}
A pairing, such that the trajectories of two nodes cross, does not minimize the total distance travelled by all the nodes.
\end{pairing}
\textbf{Proof:} Consider any two nodes $i,j$ of the initial configuration and suppose they are paired with two nodes $i',j'$ of the final configuration. There are two possible pairings, either $i$ is paired with $i'$ and $j$ is paired with $j'$ or $i$ is paired with $j'$ and $j$ is paired with $i'$, as plotted in Figure \ref{trajcross}. For the intersecting trajectories (dotted lines) the total distance
travelled by the nodes is equal to $b_1+b_2+b_3+b_4$. For the
non-intersecting trajectories (solid lines), the total distance
travelled by the nodes is equal to $a_1+a_2<b_1+b_2+b_3+b_4$, this
from the triangle inequality. Hence if the chosen pairing of nodes have intersecting trajectories of  any pair of nodes, the total distance travelled by all the nodes is not minimized. $\quad \square$

Given  initial and final configurations $\mathbf{q}_0$ and $\mathbf{q}_f$, we design a smooth path of tensegrities
$\mathbf{q}^e(t)=(\mathbf{x}^e(t),\mathbf{y}^e(t))$ using the following linear interpolation:
\begin{equation}\label{linearinterpolation}
\begin{aligned}
\mathbf{q}^e(t)&=\frac{t}{\tau}\mathbf{q}_f+(1-\frac{t}{\tau})\mathbf{q}_0,  \;\;\;\;
t\in [0, \tau]\\
\mathbf{q}^e(t)&=\mathbf{q}_f,  \;\;\;\; t>\tau.
\end{aligned}
\end{equation}
The parameter $\tau$ has units of time and allows for tuning of the tensegrity's reconfiguration speed. The greater $\tau$ is, the slower the tensegrity reconfigures itself. The linear interpolation is an attractive choice of path, limiting the distance travelled by each vehicle and hence the energy required for the reconfiguration assuming no strong steady external forces. 

\paragraph{Control law:}  We now consider the following parameterization of the system (\ref{eomfirstorder2}):
\begin{equation}\label{eomfirstorder3}
\left\{ \begin{aligned}
\dot{x}_i&=p_i^x\\
\dot{y}_i&=p_i^y\\
\dot{p}_i^x&=-\nu p_i^x-\sum_{j=1}^{N}\alpha_{ij}(t)\omega_{ij}(t)\left(1-\frac{l_{ij}(t)}{r_{ij}(t)}\right)(x_i-x_j)\\
\dot{p}_i^y&=-\nu
p_i^y-\sum_{j=1}^{N}\alpha_{ij}(t)\omega_{ij}(t)\left(1-\frac{l_{ij}(t)}{r_{ij}(t)}\right)(y_i-y_j)
\end{aligned}\right.\;\;
\end{equation}
$i=1,\ldots,N$, where the parameters $\alpha_{ij}(t),l_{ij}(t)$ are to be designed so that  the designed path $[\mathbf{q}^e](t)$ is a solution for $t \in [0,\tau]$ and for fixed $t\in[0,\tau]$ the corresponding equilibrium set $[\mathbf{z}^e](t)=([\mathbf{q}^e],\mathbf{0},\mathbf{0})(t)$ is exponentially stable. Following a similar procedure to the one developed in Section \ref{nonlinear} we pick the parameters $\alpha_{ij}(t),l_{ij}(t)$ and $\omega_{ij}(t)$ so that
\begin{equation}\label{varparam}
\begin{aligned}
\begin{pmatrix}\tilde{\Omega}(t,\mathbf{x}^e(t),\mathbf{y}^e(t)) & 0\cr 0
&\tilde{\Omega}(t,\mathbf{x}^e(t),\mathbf{y}^e(t))\end{pmatrix}\begin{pmatrix}\mathbf{x}^e(t)\cr\mathbf{y}^e(t)\end{pmatrix}=0,\;\;\;\;
 t\in[0, \tau]\\
\begin{pmatrix}\tilde{\Omega}(t,\mathbf{x}^e(\tau),\mathbf{y}^e(\tau)) & 0\cr 0
&\tilde{\Omega}(t,\mathbf{x}^e(\tau),\mathbf{y}^e(\tau))\end{pmatrix}\begin{pmatrix}\mathbf{x}^e(\tau)\cr\mathbf{y}^e(\tau)\end{pmatrix}=0,\;\;\;\; t>\tau,
\end{aligned}
\end{equation}
where the explicit dependence on the first argument $t$ in $\tilde \Omega$ comes from the dependence of $\tilde \Omega$ on the parameters $\alpha_{ij}(t), l_{ij}(t)$ and $\omega_{ij}(t)$. As a first step to solve equation (\ref{varparam}), we choose the parameters $\alpha_{ij}(t), l_{ij}(t)$ for all $i,j$ so that  $\tilde{\Omega}(\mathbf{x}^e(t),\mathbf{y}^e(t))=\Omega(t)$. This last equation is solved by setting
\begin{equation}\label{varalphal}
\begin{aligned}
\alpha_{ij}(t)&=\frac{\pi}{\arctan\omega_{ij}(t)}\\
l_{ij}(t)&=r_{ij}^e(t)\left(1-\frac{1}{\pi}\arctan\omega_{ij}(t)\right),
\end{aligned}
\end{equation}
where $r_{ij}^e(t)=\|\vec{q}_i^e(t)-\vec{q}_j^e(t)\|$. Let ${\bf r}^e(t)$ be the vector with elements $r^e_{ij}(t)$ strung together.   The parameters $\omega_{ij}(t)$ are then computed as in Section \ref{nonlinear} using the identity 
\begin{equation}\label{varomega}
\Omega(t)=\Lambda(t) D\Lambda(t)^T,
\end{equation}
where $D=\mbox{diag}\begin{pmatrix} 0 & 0 & 0 & d_4&\cdots&d_N\end{pmatrix},\; d_i>0,\; i=4,\ldots,N$ and the columns of $\Lambda(t)$ constitute a basis of orthonormal eigenvectors obtained by the Gram-Schmidt procedure on  $N$ linearly independent  vectors where  $\mathbf{x}^e(t),\mathbf{y}^e(t),\mathbf{1}$ are the first three.  This choice makes $[\mathbf{q}^e](t)$ a parameterized (by $t$) family of stable equilibrium sets for the system (\ref{eomfirstorder3}). We now use results from Lawrence and Rugh \cite{Lawrence} to show that the time-varying system (\ref{eomfirstorder3}) is well behaved.

\subsection{Boundedness and convergence}

Following the notations in \cite{Lawrence}, the following setting is considered: a system described by
\begin{equation}\label{sys}
\dot{\mathbf{z}}(t)=f(\mathbf{z}(t),\mathbf{u}(t)), \;\; \mathbf{z}(0)=\mathbf{z}_0, \;\;\;\; t\geq 0,
\end{equation}
where $\mathbf{z}(t)\in\mathbb{R}^{4N}$ is the state vector, $\mathbf{u}(t)=(\alpha_{ij}(t),l_{ij}(t),\omega_{ij}(t),r^e_{ij}(t))\in\mathbb{R}^m$ is the input vector and $f$ is the vector field given by (\ref{eomfirstorder3}). For such system, Lawrence and Rugh proved the following boundedness result \cite{Lawrence}:
\begin{lawrence} \label{lawrence} Suppose the system (\ref{sys}) satisfies
\renewcommand{\labelenumi}{H\arabic{enumi}}
\begin{enumerate}
\item $f:\mathbb{R}^{4N}\times\mathbb{R}^m\mapsto\mathbb{R}^{4N}$ is
twice differentiable,
\item there is a bounded, open set $\Gamma\subset\mathbb{R}^m$ and
a continuously differentiable function
$\mathbf{z}:\Gamma\mapsto\mathbb{R}^{4N}$ such that for each constant
input value $\mathbf{u}\in \Gamma$,
$f(\mathbf{z}(\mathbf{u}),\mathbf{u})=0$,
\item there is a $\lambda>0$ such that for each $\mathbf{u}\in
\Gamma$, the eigenvalues of $(\partial f/\partial
\mathbf{z})(\mathbf{z}(\mathbf{u}),\mathbf{u})$ have real parts no greater
than $-\lambda$.
\end{enumerate}
Then there is a $\rho^*>0$ such that given any $\rho\in[0,
\rho^*]$ and $T>0$, there exist
$\delta_1(\rho),\delta_2(\rho,T)>0$ for which the following
property holds. If a continuously differentiable input $\mathbf{u}(t)$
satisfies $\mathbf{u}(t)\in\Gamma, t\geq t_0$,
\begin{equation*}
\|\mathbf{z}_0-\mathbf{z}(\mathbf{u}(t_0))\|<\delta_1
\end{equation*}
and
\begin{equation*}
\frac{1}{T}\int_t^{t+T}\|\dot{\mathbf{u}}(\sigma)\|d\sigma<\delta_2, \;\;\;
t\geq t_0
\end{equation*} then the corresponding solution of (\ref{sys})
satisfies
\begin{equation*}
\|\mathbf{z}(t)-\mathbf{z}(\mathbf{u}(t))\|<\rho, \;\;\; t\geq t_0.
\end{equation*}
\end{lawrence}
To apply Theorem \ref{lawrence} to our system, we first show that all three conditions H1-H3 are satisfied.

The vector field (\ref{eomfirstorder3}) is proved to satisfy H1 by showing that $\Omega(t)$ and $\mathbf{r}^e(t)$ are smooth functions of $t$ and that $f(\mathbf{z}(\cdot),\mathbf{u}(\cdot))$ is a smooth function of $\Omega(\cdot)$ and $\mathbf{r}^e(\cdot)$. The path of tensegrities $\mathbf{q}^e(t)=(\mathbf{x}^e(t), \mathbf{y}^e(t))$, given by (\ref{linearinterpolation}) is a linear interpolation between $\mathbf{q}_0$ and $\mathbf{q}_f$, and hence a smooth function of time. This guarantees smoothness for $\mathbf{r}^e(t)$. The time-varying stress matrix $\Omega(t)$ is computed as 
\begin{equation*}
\Omega(t)=\Lambda(t) D\Lambda(t)^T,
\end{equation*}
where $D$ is a constant diagonal matrix and the columns of $\Lambda(t)$ constitute an orthonormal basis of $\mathbb{R}^N$ obtained through a Gram-Schmidt procedure on  $N$ linearly independent  vectors $\mathbf{x}^e(t),\mathbf{y}^e(t),\mathbf{1},$ $\mathbf{w}_4(t),\cdots,\mathbf{w}_N(t)$. The vectors obtained from a Gram-Schmidt procedure consist of linear combinations of the original set of linearly independent vectors, hence $\Lambda(t)$ and consequently $\Omega(t)$ are smooth functions of $t$. We now show that $f(\mathbf{z}(\cdot),\mathbf{u}(\cdot))$ is a smooth function of $\Omega(\cdot)$ and $\mathbf{r}^e(\cdot)$. In Section \ref{nonlinear} we noted that the choice for the parameters $\alpha_{ij}$ and $l_{ij}$ given by
\begin{align*}
\alpha_{ij}&=\frac{\pi}{\arctan\omega_{ij}}\\
l_{ij}&=r_{ij}^e\left(1-\frac{1}{\pi}\arctan\omega_{ij}\right)
\end{align*}
is such that the vector field (\ref{eomfirstorder2}) is a smooth map of $\omega_{ij}$ and $r_{ij}^e$. Hence the vector field (\ref{eomfirstorder3}) is a smooth map of $\Omega(t)$ and $\mathbf{r}^e(t)$. This concludes the proof that the system (\ref{eomfirstorder3}) satisfies H1. 

We prove that (\ref{eomfirstorder3}) satisfies  both H2 and H3 using the results from Section \ref{sec:stability}. By Theorems \ref{isolated} and \ref{stable} we know that choosing the parameters $\alpha_{ij}, l_{ij}$ and $\omega_{ij}$ satisfying (\ref{eqcond}) makes $[\mathbf{z}^e]$ an isolated exponentially stable equilibrium set of (\ref{eomfirstorder2}).  Hence choosing $\alpha_{ij}(t), l_{ij}(t)$ and $\omega_{ij}(t)$ satisfying (\ref{varparam}) for every fixed $t$ makes $[\mathbf{z}^e](t)$ an exponentially stable equilibrium set parameterized by $\Omega(t),\mathbf{r}^e(t)$, concluding the proof that (\ref{eomfirstorder3}) satisfies both H2 and H3.

Theorem \ref{lawrence} guarantees that if we start ``close enough" to $[\mathbf{q}_0^e]$ (i.e., $\|[\mathbf{z}](0)-([\mathbf{q}_0^e], \mathbf{0})  \|<\delta_1$) and the reconfiguration is not ``too fast" (i.e., $\tau$ large enough such that $\frac{1}{T}\int_t^{t+T}\|\dot{\mathbf{u}}(\sigma)\|d\sigma<\delta_2$), then the controlled reconfiguration is well behaved, i.e., $\|[\mathbf{q}^e](t)-[\mathbf{q}](t)\|$ is bounded and $[\mathbf{q}](t) \rightarrow [\mathbf{q}_f^e]$ as $t\rightarrow\infty$. We next explore the performance of this control law with a simulation example.

\subsection{Simulations}\label{simulations}

We present in this section simulation results for the reconfiguration control law given by (\ref{eomfirstorder3}). We investigate with an example the effect of the choice of $\tau$ on the shape error and on the total distance travelled by the vehicles through the reconfiguration. Shape error $e(t)$  is measured as
\begin{equation}\label{shapeerror}
e(t)=\sum_{i=1}^N(d_{iG}(t)-d_{iG}^e(t))^2,
\end{equation}
where $d_{iG}(t)$ (respectively $d_{iG}^e(t)$) is the observed (respectively planned) distance between the $i$th node and the center of mass of the formation at time $t$. 

\begin{figure}
\begin{center}
\includegraphics[width=5in]{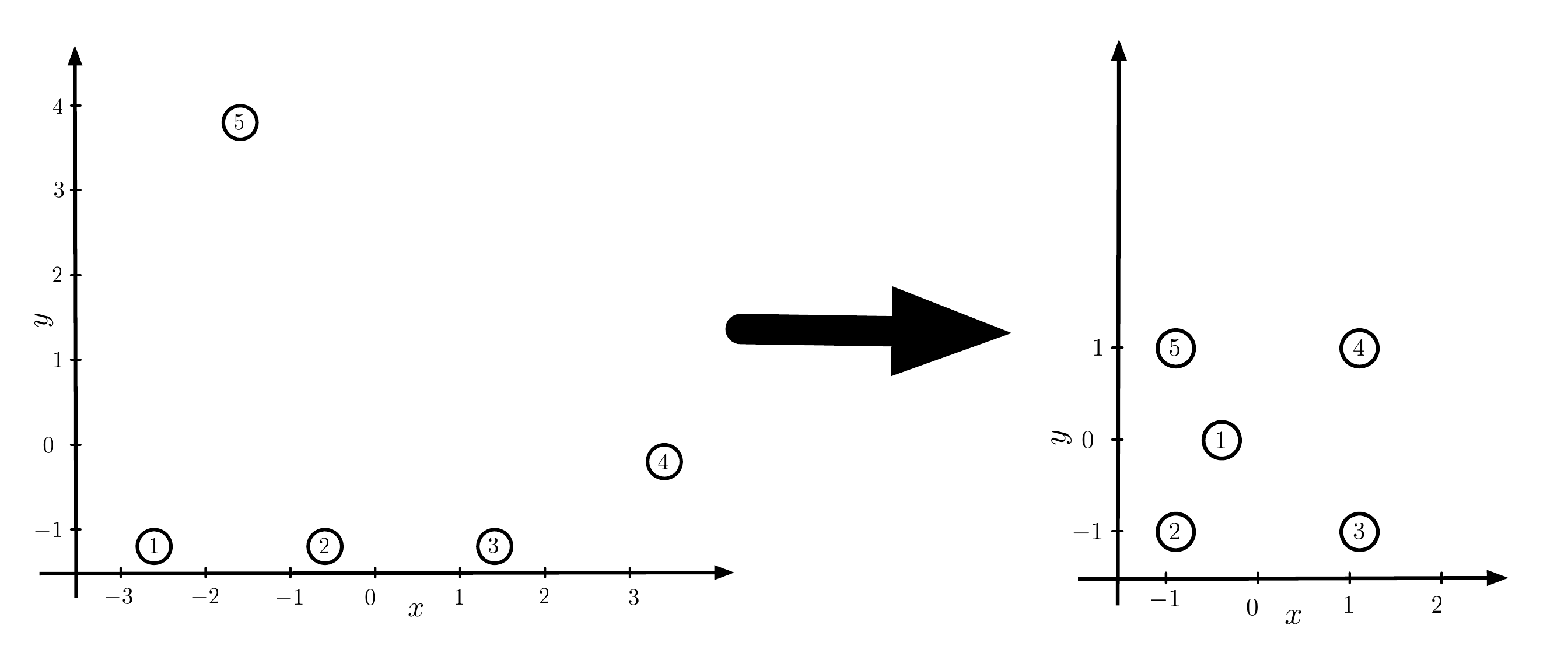}
\end{center}
\caption{Plot of the planned initial and final configurations for the studied five vehicle formation. The initial configuration is convex but not strictly convex, while the final configuration is non-convex.}
\label{initfinal}
\end{figure}

We consider as an example the reconfiguration of a five vehicle formation from initial configuration $\mathbf{q}_0$ to final configuration $\mathbf{q}_f$ given by
\begin{align*}
\mathbf{q}_0&=\begin{pmatrix} -2.6 & -.6 & 1.4 & 3.4 & -1.6\cr -1.2& -1.2& -1.2 & -.2 & 3.8\end{pmatrix}\\
\mathbf{q}_f&=\begin{pmatrix} -.4 & -.9 & 1.1 & 1.1 & -.9\cr 0& -1& -1 & 1 & 1\end{pmatrix},
\end{align*}
plotted in Figure \ref{initfinal}, where the coordinates are expressed in meters. The initial configuration is convex but not strictly convex and the final configuration is non-convex. 

\begin{figure}
\begin{center}
\includegraphics[width=3.1in]{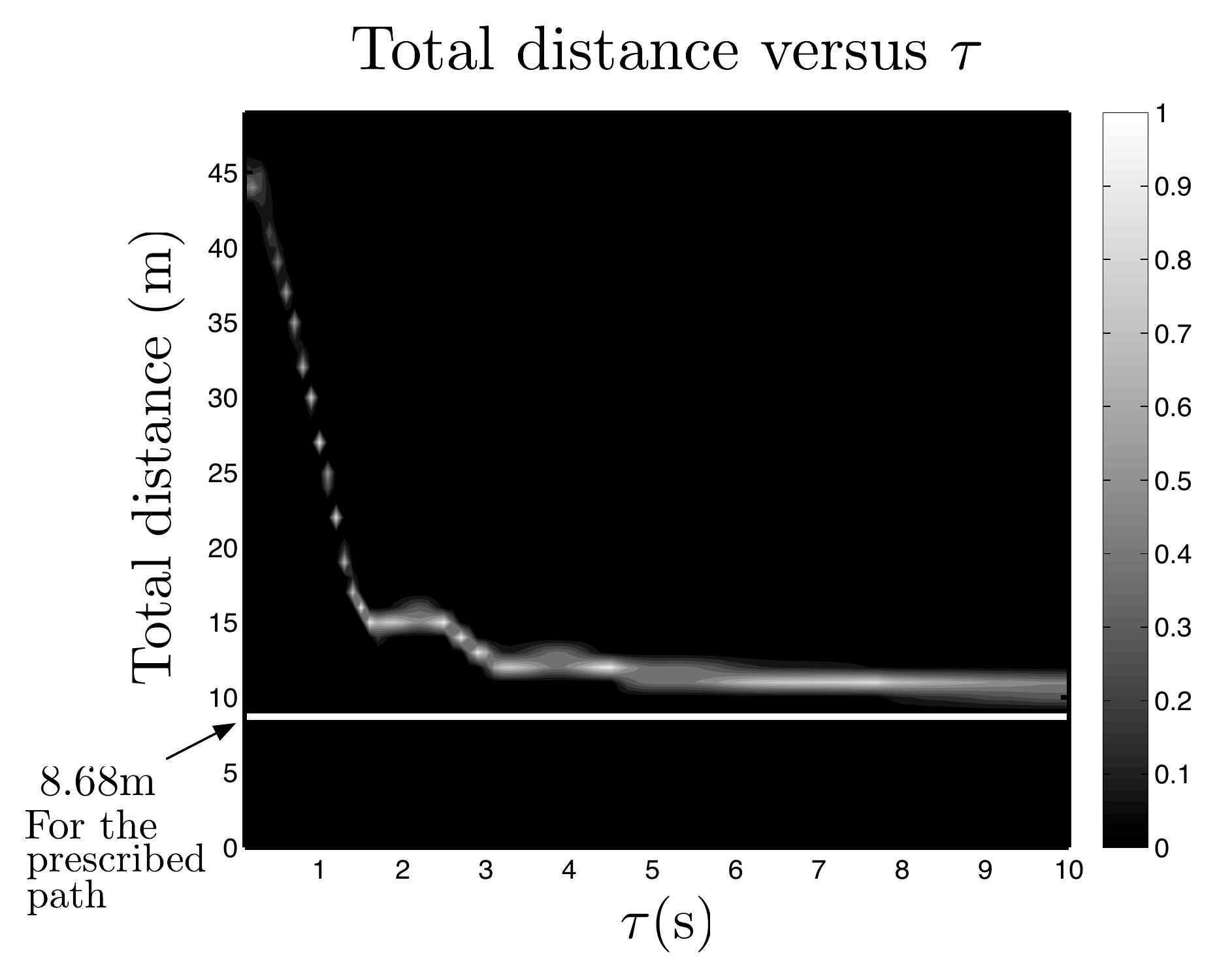}
\end{center}
\caption{Evolution of  the mean total distance travelled by the five nodes calculated over five thousand runs (each run for 500s) as a function of $\tau$, for values of $\tau$ between .1s and 10s with increment of .1. The lower bound is given by the total distance travelled for the linear designed path (solid black line which gives 8.7m).
}
\label{totaldistance}
\end{figure}

Figure \ref{totaldistance} shows the evolution of  the mean total distance travelled by the five nodes calculated over five thousand runs as a function of $\tau$, for values of $\tau$ between .1 and 10 with increment of .1.  The greater $\tau$ is (i.e., the slower the network is prescribed to reconfigure itself), the shorter the distance travelled by the five vehicles. The lower bound is given by the total distance travelled for the linear designed path (solid black line which gives 8.7m).

\begin{figure}
\begin{center}
\includegraphics[width=6.6in]{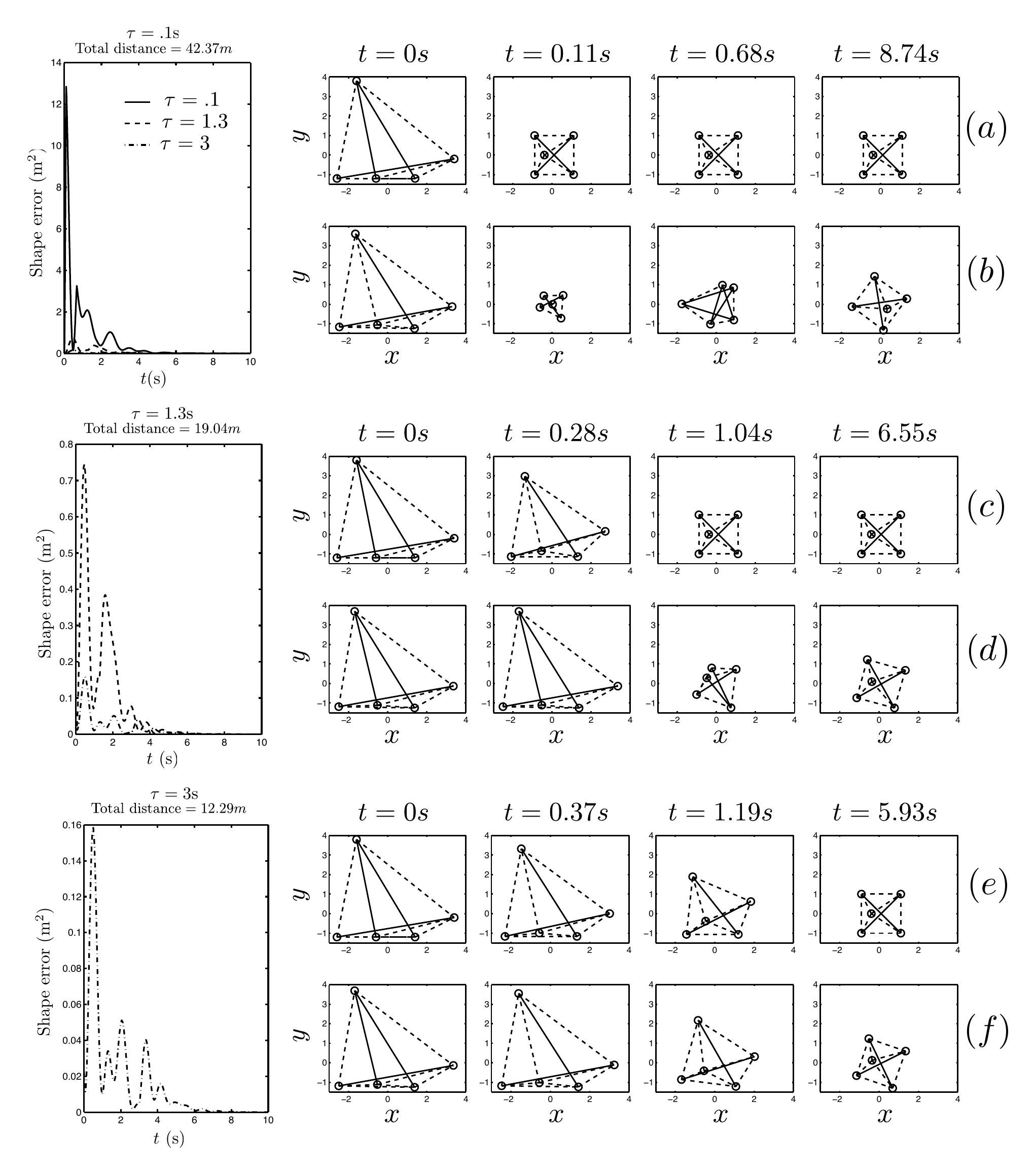}
\end{center}
\caption{For $\tau=0.1\mbox{s}, 1.3$s and $3$s, shape error is plotted as a function of time.  (a)-(c)-(e) give the snapshots of the five vehicle network along the prescribed path  respectively for $\tau=0.1\mbox{s},1.3$s and $3$s. (b)-(d)-(f)  give actual snapshot of the five vehicle network at the beginning of the reconfiguration, at the first two peak of the shape error curve and when the shape error becomes permanently smaller than $10^{-3}$m  respectively for $\tau=0.1\mbox{s},1.3$s and $3$s. }
\label{simulation}
\end{figure}

Figure \ref{simulation} shows for $\tau=.1\mbox{s} , 1.3\mbox{s}$ and $3$s a plot of the shape error as a function of time and snapshots of the five vehicle network at the beginning of the reconfiguration, at the first two peaks of the shape error curve and when the shape error becomes permanently smaller than $10^{-3}$m. We note that the shape error graphs do not all have the same scale. Looking at the graph of the shape error for the cases $\tau=.1$s and $\tau=3$s, we observe a difference of two orders of magnitude. The case where $\tau=.1$s is so fast that it is close to the case in which no intermediate points on the path are given but rather the system is required to stabilize to the final configuration given the initial configuration as initial condition, as in Section \ref{nonlinear}. Looking at the snapshots at the peak of shape error for the case $\tau=.1$s we see that these peaks can be interpreted as overshoots. Indeed the network over-shrinks at first but then over-extends whereas in the $\tau=3$s case, the structure follows more smoothly the prescribed path. Likewise, the $\tau=.1$s case takes the longest to converge and the $\tau=3$s case takes the shortest time to converge. We also note in these snapshots that, independent of the choice of $\tau$, the structure at its final configuration is rotated; this is due to changing total angular momentum of the system. It highlights the fact that we are only controlling the shape of the structure but not its position and orientation. 

\section{Final Remarks}\label{sec:finalremarks}

We have presented and proven a methodology that provides distributed control laws that stabilize a multi-vehicle formation for an arbitrary planar shape.   The results extend to changes of formation shape over a given time interval.  A key idea is to model the controlled multi-vehicle dynamics as a tensegrity structure with nodes representing vehicles and forces along struts and cables representing interconnecting vehicle control forces.   Critical to the result is the smooth map we derive from arbitrary planar shape to a parametrized tensegrity structure.  

It is a limitation, however, that given an arbitrary shape, our map often yields a tensegrity structure with number of interconnections greater than the proven lower limit.   We showed with an example how to manipulate our method to reduce the number of interconnections; however, future work will examine making such a modification systematic.    

Since our stability results are local, it is  of interest to consider proving more global results and exploring the global phase space to better leverage dynamics of tensegrity structures; for example, see \cite{NabetIFAC} for interesting, possible periodic solutions.   In \cite{Pais09} we prove global results for one-dimensional tensegrity structures and use this to prove global stability of planar and three-dimensional shapes by creating two and three-dimensional structures made up of a set of orthogonal one-dimensional tensegrity structures.   In future work we will consider extending the method of the present paper to three dimensions to provide a compelling, alternative means for controlling the shape of three-dimensional formations.

\section{Acknowledgements}
We thank Eduardo Sontag and Ming Cao for inspiring discussions on this topic.

\bibliographystyle{unsrt}
\bibliography{biblio}

\renewcommand{\thesection}{\Alph{section}}
\setcounter{section}{0}
\section{APPENDIX} \label{lemma1}
\renewcommand{\theequation}{A-\arabic{equation}}
\setcounter{equation}{0}  

\qquad In this appendix, we present the proof of Lemma~\ref{pos}. We write
$\delta^2 \tilde V(\mathbf{q}^e)$ as
\begin{equation*}
\delta^2 \tilde V(\mathbf{q}^e)=M_1+M_2:=\begin{pmatrix} \Omega& 0_N \cr
0_N & \Omega\end{pmatrix}+\begin{pmatrix} L_{\omega x}& L_{\omega
xy} \cr L_{\omega xy} & L_{\omega y}\end{pmatrix}.
\end{equation*} Recall that $\Omega$ is designed to be positive semi-definite, hence $M_1\geq 0$.
We show $M_2$ positive semi-definite. By direct computation,
\begin{equation}\label{Lwpositive}
\begin{pmatrix} \mathbf{q}_x^T & \mathbf{q}_y^T\end{pmatrix}\begin{pmatrix} L_{\omega x}& L_{\omega xy} \cr
L_{\omega xy} & L_{\omega y}\end{pmatrix}\begin{pmatrix}
\mathbf{q}_x \cr \mathbf{q}_y\end{pmatrix}=\mathbf{q}_x^T
L_{\omega x}\mathbf{q}_x+\mathbf{q}_y^T L_{\omega
y}\mathbf{q}_y+2\mathbf{q}_x^T L_{\omega xy}\mathbf{q}_y,
\end{equation} where $\mathbf{q}_x=(q_{x_1}\cdots q_{x_{N}})^T\in
\mathbb{R}^N,$ and $\mathbf{q}_y=(q_{y_1}\cdots q_{y_{N}})^T\in
\mathbb{R}^N$. Each term of the sum can be rewritten as
\begin{align*}
\mathbf{q}_x^T L_{\omega
x}\mathbf{q}_x&=\sum_{i=1}^{N}\sum_{j=1,j\neq
i}^{N}q_{x_i}\frac{\alpha_{ij}\omega_{ij}(x_i-x_j)^2l_{ij}}{r_{ij}^3}(q_{x_i}-q_{x_j})\\
&=\sum_{i<j}\frac{\alpha_{ij}\omega_{ij}(x_i-x_j)^2l_{ij}}{r_{ij}^3}(q_{x_i}-q_{x_j})^2,\\
\mathbf{q}_y^T L_{\omega
y}\mathbf{q}_y&=\sum_{i=1}^{N}\sum_{j=1,j\neq
i}^{N}q_{y_i}\frac{\alpha_{ij}\omega_{ij}(y_i-y_j)^2l_{ij}}{r_{ij}^3}(q_{y_i}-q_{y_j})\\
&=\sum_{i<j}\frac{\alpha_{ij}\omega_{ij}(y_i-y_j)^2l_{ij}}{r_{ij}^3}(q_{y_i}-q_{y_j})^2,\\
\mathbf{q}_x^T L_{\omega
xy}\mathbf{q}_y&=\sum_{i=1}^{N}\sum_{j=1,j\neq
i}^{N}q_{x_i}\frac{\alpha_{ij}\omega_{ij}(x_i-x_j)(y_i-y_j)l_{ij}}{r_{ij}^3}(q_{y_i}-q_{y_j})\\
&=\sum_{i<j}\frac{\alpha_{ij}\omega_{ij}(x_i-x_j)(y_i-y_j)l_{ij}}{r_{ij}^3}(q_{y_i}-q_{y_j})(q_{x_i}-q_{x_j}).
\end{align*}
(\ref{Lwpositive}) can now be factored as
\begin{equation}\label{Lwpositive2}
\mathbf{q}_x^T L_{\omega x}\mathbf{q}_x+\mathbf{q}_y^T L_{\omega
y}\mathbf{q}_y+2\mathbf{q}_x^T L_{\omega
xy}\mathbf{q}_y=\sum_{i<j}\frac{\alpha_{ij}\omega_{ij}l_{ij}}{r_{ij}^3}\Big((y_i-y_j)(q_{y_i}-q_{y_j})+(x_i-x_j)(q_{x_i}-q_{x_j})\Big)^2\geq
0.
\end{equation} This concludes the proof that $M_2\geq 0$, and hence $\delta^2V(\mathbf{q}^e)\geq 0$ $\quad \square$

\section{APPENDIX}\label{lemma2}
\renewcommand{\theequation}{B-\arabic{equation}}
\setcounter{equation}{0}  

\qquad In this appendix, we present the proof of Lemma~\ref{ker}. By Lemma~\ref{pos}, $M_1$ and $M_2$ are symmetric, positive semi-definite
matrices. Hence,
\begin{equation}\label{kernelcond}
\mathbf{q}\in
\mbox{ker}(\delta^2V(\mathbf{q}^e))\Longleftrightarrow
\mathbf{q}\in \mbox{ker}(M_1)\;\;  \mbox{and}\;\; \mathbf{q}\in
\mbox{ker}(M_2).
\end{equation} By design,
$\mbox{ker}(M_1)=\{\mathbf{w}_1,\mathbf{w}_2,\mathbf{w}_3\}$; by
direct computation, we check
$\{\mathbf{w}_1,\mathbf{w}_2,\mathbf{w}_3\}\in \mbox{ker}(M_2)$.
Using (\ref{kernelcond}), we conclude that the kernel of
$\delta^2V(\mathbf{q}^e)$ is exactly spanned by
$\{\mathbf{w}_1,\mathbf{w}_2,\mathbf{w}_3\}$. $\quad \square$

\end{document}